\documentclass[11pt]{article}
\usepackage{amsfonts}
\usepackage{graphics}
\usepackage{indentfirst}
\usepackage{cite}
\usepackage{latexsym}
\usepackage[colorlinks, linkcolor=red]{hyperref}
\usepackage{color, xcolor}
\usepackage[paper=a4paper, left=1.6cm, right=1.6cm, top=1.8cm, bottom=1.6cm, headheight=5.5pt, footskip=0.8cm, footnotesep=0.8cm, centering, includefoot]{geometry}
\hypersetup{urlcolor=red, citecolor=red}
\usepackage{amsmath}
\allowdisplaybreaks
\usepackage{amssymb}
\usepackage[dvips]{epsfig}
\usepackage{amscd}

\newtheorem{theorem}{Theorem}[section]
\newtheorem{remark}{Remark}[section]

\newtheorem{lemma}{Lemma}[section]

\newtheorem{proposition}{Proposition}[section]

\DeclareMathOperator{\divv}{div}
\DeclareMathOperator{\curl}{curl}

\makeatletter
\@addtoreset{equation}{section}
\makeatother
\makeatletter
\@addtoreset{equation}{section}
\makeatother

\title{Local well-posedness to the 2D Cauchy problem of full compressible magnetohydrodynamic equations with vacuum at infinity
\thanks{This research was partially supported by National Natural Science Foundation of China (Nos. 11901474, 12071359), Exceptional Young Talents Project of Chongqing Talent (No. cstc2021ycjh-bgzxm0153), and the Innovation Support Program for Chongqing Overseas Returnees (No. cx2020082).}
}

\date{}
\author{ Hong Chen\thanks{School of Mathematics and Statistics, Southwest University, Chongqing 400715,
People's Republic of China ({\tt chenhong123@office365.swu.edu.cn}). }
\quad  Xin Zhong\thanks{Corresponding author. School of Mathematics and Statistics, Southwest University, Chongqing 400715,
People's Republic of China ({\tt xzhong1014@amss.ac.cn}).
}
}

\begin{document}
\maketitle

\begin{abstract}
This paper concerns the Cauchy problem of two-dimensional (2D) full compressible magnetohydrodynamic (MHD) equations in the whole plane $\mathbb{R}^2$ with zero density at infinity. By spatial weighted energy method, we derive the local existence and uniqueness of strong solutions provided that the initial density and the initial magnetic field decay not too slowly at infinity. Note that the initial temperature does not need to decay slowly at infinity. In particular, vacuum states at both the interior domain and the far field are allowed.
\end{abstract}



\section{Introduction}

The motion of a $n$-dimensional compressible viscous, heat-conductive, magnetohydrodynamic flows is governed by the following full compressible MHD
equations (see \cite[Chapter 3]{TT2012}):
\begin{align}\label{mhd}
\begin{cases}
\rho_t + \divv(\rho u) = 0,\\
(\rho u)_t+\divv(\rho u\otimes u)-\mu\Delta u
-(\mu+\lambda)\nabla\divv u+\nabla P=H\cdot\nabla H-\frac{1}{2}\nabla |H|^2,\\
c_v[(\rho \theta)_{t}+\divv(\rho u\theta)]-\kappa\Delta\theta+P\divv u =\frac{\mu}{2}|\nabla u+(\nabla u)^{tr}|^{2}+\lambda(\divv u)^2
+\nu|\curl H|^{2},\\
H_t-H\cdot\nabla u+u\cdot\nabla H+H\divv u=\nu\Delta H,\\
\divv H=0.
\end{cases}
\end{align}
Here $t\geq0$ is the time, $x\in \Omega \subset \Bbb R^n$ is the spatial coordinate, $\rho=\rho(x,t)$, $u=(u^1,u^2,\cdots,u^n)(x,t)$, $\theta=\theta(x,t)$, and $H=(H^1,H^2,\cdots,H^n)(x,t)$ denote the density, velocity, absolutely temperature, and the magnetic field, respectively. The pressure $P$ is given by
\begin{align*}
P(\rho)=R\rho\theta,
\end{align*}
where $R$ is a positive constant. The constant viscosity coefficients $\mu$ and $\lambda$ satisfy the physical restrictions
 \begin{align*}
 \mu>0,\ 2\mu+n\lambda\geq0.
 \end{align*}
The positive constants $c_v$ and $\kappa$ are the heat capacity and the ratio of the heat conductivity coefficient over the heat capacity, respectively, while $\nu>0$ is the magnetic diffusivity.

Let $\Omega=\mathbb{R}^2$ and we consider the Cauchy problem of \eqref{mhd} with the initial condition
\begin{align}\label{a2}
(\rho, \rho u, \rho\theta, H)(x, 0)=(\rho_0, \rho_0u_0,\rho_0\theta_0 , H_0)(x), \ \ x\in\Bbb R^2,
\end{align}
and the far field behavior
\begin{align}\label{n4}
(\rho, u, \theta, H)(x, t)\rightarrow (0, 0, 0, 0)\ \ {\rm as}\ |x|\rightarrow \infty, \ t>0.
\end{align}

The aim of the present paper is to study the existence and uniqueness of strong solutions to the 2D Cauchy problem \eqref{mhd}--\eqref{n4}.
The initial density is allowed to vanish and the spatial measure of the set of vacuum can be arbitrarily large, in particular, \textit{the initial density can even have compact support}. This study is motivated strongly by the local well-posedness of strong solutions for full compressible Navier-Stokes equations \cite{L2021} which we improve at some technical
points. The well-posedness of solutions for full compressible Navier-Stokes equations is a classical problem in the mathematical theory of fluid dynamics and has been studied extensively (see, for example, \cite{H97,L21,cho1,LWX19,DM08,HL18,CTWZ20,YZ20,XZ21,
MN83,WZ17,JL19,LZ21,L19,LX20,LX22}). Yet the well-posedness theory of multi-dimensional problems becomes much more complicated and challenging in the presence of vacuum (that is, the flow density is zero) despite its fundamental importance both physically and theoretically in understanding the
behaviors of the solutions of compressible viscous flows. Substantial difficulties arise due to the strong degeneracy of the hyperbolic-parabolic system in the case where a vacuum state appears (see \cite{L21,L2021,cho1,HL18,YZ20,WZ17,JL19,LZ21}).

When the heat conduction can be neglected and the compressible viscous fluids
are isentropic, the full compressible MHD equations \eqref{mhd} can be reduced to
the following system
\begin{align}\label{1.4}
\begin{cases}
\rho_t + \divv(\rho u) = 0,\\
(\rho u)_t+\divv(\rho u\otimes u)-\mu\Delta u
-(\mu+\lambda)\nabla\divv u+\nabla P=H\cdot\nabla H-\frac{1}{2}\nabla |H|^2,\\
H_t-H\cdot\nabla u+u\cdot\nabla H+H\divv u=\nu\Delta H,\\
\divv H=0,
\end{cases}
\end{align}
where the equation of state satisfies
\begin{align*}
P=A\rho^\gamma,\ A>0, \ \gamma>1.
\end{align*}
The well-posedness problem for the multi-dimensional compressible MHD equations in the presence of vacuum is important both in theory and in applications. Due to the strong coupling between the fluid motion and the magnetic field, the mathematical study of a compressible MHD system is rather complicated. Now, we briefly recall some existence results concerning the compressible MHD equations \eqref{1.4} with vacuum. In terms of Lions-Feireisl compactness framework for compressible Navier-Stokes equations \cite{L98,F04}, Hu and Wang \cite{HW2010} proved global weak solutions to the 3D initial boundary value problem with finite energy for $\gamma>\frac32$. Non-uniqueness of global-in-time weak solutions for an inviscid fluid in two dimensions was investigated by Feireisl and Li \cite{FL2020}. On the other hand, under the following compatibility condition
\begin{align*}
-\mu\Delta u_0-(\mu+\lambda)\nabla\divv u_0+\nabla P(\rho_0)-H_0\cdot\nabla H_0+\frac{1}{2}\nabla |H_0|^2=\sqrt{\rho_0}g
\end{align*}
for some $g\in L^2(\mathbb{R}^3)$,
Li-Xu-Zhang \cite{LXZ2013} established the global existence and uniqueness of classical solutions to the 3D Cauchy problem with smooth initial data which are of small energy but possibly large oscillations and vacuum states at both the interior domain and the far field. The study is motivated by a similar result for global classical solutions of the isentropic compressible Navier-Stokes equations obtained by Huang-Li-Xin \cite{HLX12}. Later, Hong-Hou-Peng-Zhu \cite{HHPZ2017} improved such result by allowing the initial energy large as long as $\gamma$ is close to 1 and $\nu$ is suitably large. Meanwhile, by the spatial weighted energy method, L{\"u}-Shi-Xu \cite{LSX2016} showed the global existence and uniqueness of strong solutions to the 2D Cauchy problem provided that the smooth initial data are of small total energy.

Let's turn our attention to the full compressible MHD system \eqref{mhd} with vacuum.
Compared with the isentropic case \eqref{1.4}, the additional difficulty
for studying the well-posedness of \eqref{mhd}
is that the basic energy inequality does not provide any useful dissipation
estimates on $u$ and $H$. Using the entropy method, Ducomet and Feireisl \cite{BE2006} studied the global existence of weak solutions by introducing the entropy equation rather than the thermal equation \eqref{mhd}$_3$ under the assumption
that the viscosity coefficients depend on the temperature and the magnetic field. Meanwhile, Hu and Wang \cite{XD2008} proved global weak solutions with finite energy and temperature-dependent heat conductivity.
Moreover, Li and Sun \cite{LS2021} showed global-in-time weak solutions with large initial data for the 2D non-resistive case (that is, $\nu=0$ in \eqref{mhd}). Yet the uniqueness and regularity of such weak solutions is still open. On the other hand, under some compatibility condition for the initial data, Fan and Yu \cite{JW2009} derived the local existence theory of strong solutions to the 3D problem. Later on, Huang and Li \cite{HL13} showed that if $T^*$ is the finite blow up time for the solution obtained in \cite{JW2009}, then the following Serrin type criterion holds true
\begin{equation}\label{1.5}
\lim_{T\rightarrow T^{*}}\left(\|\rho\|_{L^\infty(0,T;L^\infty)}+\|u\|_{L^s(0,T;L^r)}\right)
=\infty,\ \text{for}\ \ \frac2s+\frac3r\leq1,\ s>1,\ 3<r\leq\infty.
\end{equation}
Recently, applying \eqref{1.5} and delicate energy estimates, Liu and Zhong \cite{LZ2020} extended the local strong solutions for the 3D Cauchy problem to be a global one provided that $\|\rho_0\|_{L^\infty}+\|H_0\|_{L^3}$ is suitably small and the viscosity coefficients satisfy $3\mu>\lambda$. Meanwhile, by similar strategies as those in \cite{LZ2020}, Hou-Jiang-Peng \cite{HJP22} established global-in-time existence of strong solutions under the assumption that $\|\rho_0\|_{L^1}+\|H_0\|_{L^2}$ is suitably small. Very recently, Liu and Zhong \cite{LZ22} improved such global existence results and showed the global well-posedness of strong solutions as long as the initial energy is small enough. Moreover, they also obtained
algebraic decay estimates of the solution. However, it is still open even for the local existence of strong solutions to the 2D Cauchy
problem \eqref{mhd}--\eqref{n4}. In fact, this is the main goal of this paper.

Our main result can be stated as follows.
\begin{theorem}\label{t1}
Let $\eta_0 $ be a positive constant and
\begin{align}\label{2.07}
\bar x\triangleq(e+|x|^2)^{\frac12}\ln^{1+\eta_0} (e+|x|^2).
\end{align}
For constants $q>2$ and $a>1$, assume that the initial data
$(\rho_0\geq0, u_0, \theta_0\geq0, H_0)$ satisfies
\begin{align}\label{1.9}
\begin{cases}
\rho_0\bar x^a \in L^1 \cap H^1\cap W^{1,q},
\ H_0\bar x^{\frac{a}{2}}\in H^1,\ \divv H_0=0, \\
\big(\sqrt{\rho_0}u_0,\sqrt{\rho_0}\theta_0\big)\in L^2,\
\big(\nabla u_0, \nabla \theta_0, \nabla H_0\big)\in H^1,
\end{cases}
\end{align}
and the compatibility condition
\begin{align}\label{c.c}
\begin{cases}
-\mu\Delta u_0-(\mu+\lambda)\nabla\divv u_0+\nabla(R\rho_0\theta_0)-H_0\cdot\nabla H_0+\frac{1}{2}\nabla |H_0|^2=\sqrt{\rho_0}g_1, \\
\kappa\Delta\theta_0+\frac{\mu}{2}|\nabla u_0+(\nabla u_0)^{tr}|^2
+\lambda(\divv u_0)^2+\nu(\curl H_0)^2-R\rho_0\theta_0\divv u_0=\sqrt{\rho_0}g_2,
\end{cases}
\end{align}
for some $g_1,g_2\in L^2(\mathbb{R}^2)$, where $\curl H\triangleq\partial_1H^2-\partial_2H^1$.
Then there exists a small time $T_0>0$ such that the problem  \eqref{mhd}--\eqref{n4} has a unique strong solution $(\rho\geq0,u,\theta\geq0,H)$ on $\mathbb{R}^2\times (0,T_0]$ satisfying
\begin{align}\label{1.10}
\begin{cases}
\rho\in C([0,T_0];L^1 \cap H^1\cap W^{1,q}),\\
\rho\bar x^a\in L^\infty(0,T_0;L^1\cap H^1\cap W^{1,q}),\\
\sqrt{\rho}u,\sqrt{\rho}\theta, \sqrt{\rho}u_t,
\sqrt{\rho}\theta_t\in L^\infty(0,T_0; L^2), \\
\nabla u,\ \nabla\theta,\ H\bar x^{\frac{a}{2}}\in\,L^\infty(0,T_0; H^1), \\
H, \nabla H,  H_t, \nabla^2 H\in L^\infty(0,T_0;L^2), \\
\nabla u, \nabla\theta\in L^2(0,T_0; W^{1,q})\cap  L^{\frac{q+1}{q}}(0,T_0; W^{1,q}), \\
H_t, \nabla H\bar{x}^{\frac{a}{2}}\in L^2(0, T_0; H^1),\\
\sqrt{\rho}u_t, \sqrt{\rho}\theta_t, \nabla u_t,
\nabla\theta_t\in L^2(0, T_0;L^2),\\
\end{cases}
\end{align}
and
\begin{align}\label{l1.2}
\inf\limits_{0\le t\le T_0}\int_{B_{N}}\rho(x,t)dx\ge \frac14\int_{\mathbb{R}^2} \rho_0(x)dx,
\end{align}
for some constant $N >0$ and $B_{N }\triangleq\left.\left\{x\in\mathbb{R}^2\right|
\,|x|<N \right\}$.
\end{theorem}

\begin{remark}
The compatibility condition \eqref{c.c} is used for obtaining the $L^\infty(0, T; L^2)$-norms of $\sqrt{\rho}u_t$ and $\sqrt{\rho}\theta_t$, which is crucial in dealing with the $L^\infty(0, T; L^2)$-norm of the gradient of the temperature.
This is very different from the case of isentropic flows \cite{LB2015}, where the authors showed the local well-posedness of strong solutions without using such compatibility condition via time weighted techniques. It is natural to investigate the local existence  of strong solutions to \eqref{mhd}--\eqref{n4} without using \eqref{c.c}, and this will be left for the future study.
\end{remark}

\begin{remark}
We point out that the finiteness of initial mass is a crucial assumption, which is used to obtain a Hardy type inequality for the velocity and temperature (see \eqref{3.v2}).
It is still unknown whether this condition is necessary for the local existence of strong solutions to the 2D Cauchy problem with vacuum at infinity.
\end{remark}

We now make some comments on the key analysis of this paper.
It should be noted that the crucial arguments for 3D case \cite{JW2009} cannot be applied here. Precisely speaking, one can get $u\in L^6(\mathbb{R}^3)$ provided that $\nabla u\in L^2(\mathbb{R}^3)$ and $\lim\limits_{|x|\rightarrow\infty} u=0$. However, due to the presence of vacuum at infinity and the criticality of Sobolev's inequality in $\mathbb{R}^2$, it seems impossible to control $u\in L^p(\mathbb{R}^2)$ for any $p>1$ in just terms of $\nabla u\in L^2(\mathbb{R}^2)$ and $\sqrt{\rho} u\in L^2(\mathbb{R}^2)$.
Moreover, compared with the case of 2D isentropic Cauchy problem \cite{LB2015}, some new difficulties arise due to the appearance of energy equation \eqref{mhd}$_3$ as well as the coupling of the velocity with the temperature. Indeed, if we multiply \eqref{mhd}$_3$ by $\theta$, we get after integration by parts that
\begin{align}\label{1.12}
& \frac{c_v}{2}\frac{d}{dt}\int\rho\theta^2dx+\kappa\int|\nabla\theta|^2dx \notag \\
& = -R\int \rho\theta^2\divv udx+\int\Big[\frac{\mu}{2}|\nabla u+(\nabla u)^{tr}|^{2}+\lambda(\divv u)^2+\nu(\curl H)^{2}\Big]\theta dx.
\end{align}
Since the $L^p(\mathbb{R}^2)$-norm of $\theta$ and spatial weighted estimates
on the gradients of the velocity and the magnetic field are unavailable, it is hard to control the term on the right hand side of \eqref{1.12} directly.
To this end, motivated by \cite{L2021}, we obtain a spatial weight estimate (see \eqref{2.10}) on the quadratic nonlinearity
$\frac{\mu}{2}|\nabla u+(\nabla u)^{tr}|^{2}+\lambda(\divv u)^2+\nu(\curl H)^{2}$,
which plays a crucial role in tackling the \textit{a priori} estimates of the temperature (see Lemma \ref{lemma 3.8}). Furthermore, it should be emphasized that a Hardy type inequality (see \eqref{3.i2}) and Gagliardo-Nirenberg inequality (see \eqref{2.4} and \eqref{2.5}) are mathematically useful for the analysis.

The rest of the paper is organized as follows. In Section \ref{sec2}, we collect some
elementary facts and inequalities which will be needed in later analysis. Sections \ref{sec3} is devoted to the a priori estimates which are needed to obtain  the local existence and uniqueness of strong solutions. Finally, the main result Theorem \ref{t1} is proved in Section \ref{sec4}.

\section{Preliminaries}\label{sec2}

In this section, we recall some known facts and elementary
inequalities which will be used later.

First of all, if the initial density is strictly away from vacuum, the following local existence theorem on bounded balls can be shown by similar arguments as those in \cite{JW2009,K1983}.
\begin{lemma}\label{th0}
For $R>0$ and $B_R=\{x\in\mathbb{R}^2 ||x|<R\}$, assume that $(\rho_0,u_0,\theta_0,H_0)$ satisfies
\begin{align}\label{2.1}
\begin{cases}
 \inf\limits_{x\in B_R}\rho_0(x) >0,\,\,\rho_0\in H^2(B_R),\,\,(u_0,\theta_0,H_0)\in H_0^1(B_R)\cap H^2(B_R),\\
(\rho,u,\theta,H)(x,t=0)=(\rho_0,u_0,\theta_0,H_0)(x),\,\,\divv H_0=0,\,\,x\in  B_R\\
(u,\theta,H)=(0,0,0),\ x\in \partial B_R,\ t>0.\\
\end{cases}
\end{align}
Then there exists a small time $T_R>0$ and a unique strong solution $(\rho>0, u, \theta,H)$ to the following initial-boundary-value problem \eqref{mhd} and \eqref{2.1} on
$B_R\times(0,T_R]$ such that
\begin{align}\label{2.2}
\begin{cases}
\rho\in C([0,T_R];H^2),\\
(u,\theta, H)\in C([0,T_R];H_0^1\cap H^2)\cap L^2(0,T_R;H^3),\\
(u_t,\theta_t, H_t)\in L^\infty(0,T_R;L^2)\cap L^2(0,T_R;H^1),
\end{cases}
\end{align}
where we denote $H^k=H^k(B_R)$ for positive integer $k$.
\end{lemma}

Next, the following well-known Gagliardo-Nirenberg inequality
(see \cite[Chapter II]{LSU1968}) will be used in the next section frequently.
\begin{lemma}
For $f\in H^1(B_R)$ and $g\in L^r(B_R)\cap W^{1, q}(B_R)$ with $r\in(1,\infty)$ and $q\in(2,\infty)$, there exists a positive constant $C$ independent of $R$ such that
\begin{align}\label{2.4}
&\|f\|^p_{L^p}\leq C\|f\|^{2}_{L^2}\|f\|^{p-2}_{H^1}, \ \forall p\in[2,\infty),\\ \label{2.5}
&\|g\|_{L^\infty}\leq
C\|g\|_{L^r}+C\|g\|^{\frac{r(q-2)}{2q+r(q-2)}}_{L^r} \|\nabla g\|^{\frac{2q}{2q+r(q-2)}}_{L^q}.
\end{align}
\end{lemma}

Next, for $\Omega=\mathbb{R}^2$ or $\Omega=B_R$, the following weighted $L^m$-bounds for elements of the Hilbert space $\widetilde{D}^{1,2}(\Omega)\triangleq\{v\in H^1_{\rm loc}(\Omega)|\nabla v\in L^2(\Omega)\}$ can be found in \cite[Theorem B.1]{L1996}.
\begin{lemma}\label{1leo}
For $m\in [2,\infty)$ and $\theta\in (1+\frac{m}{2},\infty),$ there exists a positive constant $C$ such that for either $\Omega=\mathbb{R}^2$ or $\Omega=B_R$ with $R\ge 1$ and for any $v\in \widetilde{D}^{1,2}(\Omega)$,
\begin{align}\label{3h}
\left(\int_{\Omega} \frac{|v|^m}{(e+|x|^2)\ln^{\theta}(e+|x|^2)}dx  \right)^{\frac1m}\le C\|v\|_{L^2(B_1)}+C\|\nabla v\|_{L^2 (\Omega)}.
\end{align}
\end{lemma}

A useful consequence of Lemma \ref{1leo} is the following crucial weighted  bounds (see \cite[Lemma 2.4]{L2021}).
\begin{lemma}\label{lemma2.3}
Let $\bar x$ and $\eta_0$ be as in \eqref{2.07} and $\Omega$ be as in Lemma \ref{1leo}. Assume that $\rho\in L^1(\Omega)\cap L^\infty(\Omega)$ is a non-negative function such that
\begin{align}\label{2.i2}
\int_{B_{N_1}}\rho dx\ge M_1,
\end{align}
for positive constants $M_1$ and $N_1\ge 1$ with $B_{N_1}\subset\Omega.$ Then, for $\varepsilon,\eta>0$, there is a positive constant $C$ depending only on $\varepsilon,\eta, M_1,
N_1$ and $\eta_0$ such that, for $v\in \widetilde{D}^{1,2}(\Omega)$ with $\sqrt{\rho}v\in L^2(\Omega)$,
\begin{align}\label{3.i2}
\|v\bar x^{-\eta}\|_{L^{\frac{2+\varepsilon}{\tilde{\eta}}}(\Omega)}
&\le C\|\sqrt{\rho}v\|_{L^2(\Omega)}+C\big(1+\|\rho\|_{L^\infty(\Omega)}\big)\|\nabla v\|_{L^2(\Omega)}
\end{align}
with $\tilde{\eta}=\min\{1,\eta\}$.
\end{lemma}

Next, the following $L^p$-bound for the Lam{\'e} system and elliptic equations, whose proof is similar to that of \cite[Lemma 12]{CCK04}, is a direct result of the combination of the well-known elliptic theory \cite{adn2} and a standard scaling procedure.
\begin{lemma} \label{lemma2.4}
Let $v,w\in W_0^{1,p}(B_R) \ (p>1)$ be weak solutions of
\begin{align*}
\begin{cases}
\mu\Delta v+(\mu+\lambda)\nabla\divv v=F, & x\in B_R,\\
v=0, & x\in\partial B_R,
\end{cases}
\end{align*}
and
\begin{align*}
\begin{cases}
\Delta w=G, & x\in B_R,\\
w=0, & x\in\partial B_R,
\end{cases}
\end{align*}
respectively. If $F,G \in L^p(B_R)$, then there exists a positive constant $C$ independent of $R$ such that
\begin{align}\label{lp}
\|\nabla^2v\|_{L^p(B_R)}\le C\|F\|_{L^p(B_R)},\ \|\nabla^2w\|_{L^p(B_R)}\le C\|G\|_{L^p(B_R)}.
\end{align}
\end{lemma}

Finally, by the same arguments as those in \cite[Lemma 3.1]{L2021}, we have the following spatial weighted estimate on the solution.
\begin{lemma}\label{lemma2.5}
Let $(\rho,u,\theta,H)$ be the solution to the problem \eqref{mhd} and \eqref{2.1}, then it holds that, for $b_1>0$,
\begin{align}\label{2.10}
& \int_{B_R} \Big[\frac{\mu}{2}|\nabla u+(\nabla u)^{tr}|^{2}+\lambda(\divv u)^2+\nu(\curl H)^{2}\Big]|x|^{b_1}dx \notag \\
& \leq \int_{B_R}|c_v(\rho\theta_t+\rho u\cdot\nabla\theta)+P\divv u||x|^{b_1}dx.
\end{align}
\end{lemma}

\section{\textit{A priori} estimates}\label{sec3}

In this section, for $r\in [1,\infty]$ and $k\ge0$, we write
\begin{align*}
\int \cdot dx=\int_{B_R}\cdot dx, \ L^r=L^r(B_R),\ W^{k,r}=W^{k,r}(B_R),\ H^k=W^{k,2}.
\end{align*}
Moreover, for $R>4N_0\ge4$ with $N_0$ fixed, assume that $(\rho_0,u_0,\theta_0,H_0)$ satisfies, in addition to \eqref{2.1}, that
\begin{align}\label{w1}
\frac12\le \int_{B_{N_0}}\rho_0(x)dx \le \int_{B_R}\rho_0(x)dx \le \frac32.
\end{align}
Lemma \ref{th0} thus yields that there exists some $T_R>0$ such that the  initial-boundary-value problem \eqref{mhd} and \eqref{2.1} has a unique solution $(\rho>0,u,\theta,H)$ on $B_R\times(0,T_R]$ satisfying \eqref{2.2}.

Let $\bar{x}, \eta_0, a$ and $q$ be as in Theorem \ref{t1}, the main aim of this section is to derive the following key \textit{a priori} estimate on $\psi $ defined by
\begin{align}\label{3.2}
\psi(t)\triangleq & 1+\|\sqrt{\rho}u\|_{L^2}+\|\sqrt{\rho}\theta\|_{L^2}
+\|\sqrt{\rho}u_t\|_{L^2}+\|\sqrt{\rho}\theta_{t}\|_{L^2}\notag \\
\quad & +\|\nabla u\|_{H^1}+\|\nabla \theta\|_{H^1}+\|H\|_{H^2} +\|\bar{x}^{\frac{a}{2}}H\|_{H^1}+\|\rho\bar{x}^{a}\|_{L^{1}\cap H^{1}\cap W^{1,q}}.
\end{align}

\begin{proposition} \label{pro}
Assume that $(\rho_0,u_0,\theta_0,H_0)$ satisfies \eqref{c.c} and \eqref{w1}. Let $(\rho,u,\theta,H)$ be the solution to the initial-boundary-value problem \eqref{mhd} and \eqref{2.1} on $B_R\times (0,T_R]$ obtained by Lemma \ref{th0}. Then there exist positive constants $T_0$ and $M$ both depending only on $\mu, \lambda, \nu, c_v, \kappa, q$, $a$, $\eta_0$, $N_0,$ and $E_0$ such that
\begin{align}\label{o1}
& \sup\limits_{0\le t\le T_0}\big(\psi(t)+\|H_t\|_{L^2}\big)
+\int_0^{T_0}
\big(\|\sqrt{\rho}u_t\|_{L^2}^2+\|\sqrt{\rho}\theta_t\|_{L^2}^2+\|H_t\|_{L^2}^2+\|\nabla H\bar x^{\frac{a}{2}}\|_{L^2}^2\big) dt \notag \\
& \quad +\int_0^{T_0} \Big(\|\nabla^2u\|_{L^q}^{\frac{q+1}{q}}
+\|\nabla^2\theta\|_{L^q}^{\frac{q+1}{q}}
+\|\nabla^2u\|_{L^q}^{2}+\|\nabla^2\theta\|_{L^q}^{2}\Big)dt \notag \\
& \quad +\int_0^{T_0}\big(\|\nabla u_t\|_{L^2}^{2}
+\|\nabla\theta_t\|_{L^2}^{2}+\|\nabla H_t\|_{L^2}^{2}
+\|\nabla^2H\bar x^{\frac{a}{2}}\|_{L^2}^2\big)dt
\le M,
\end{align}
where
\begin{align*}
E_0\triangleq& \|\sqrt{\rho_0}u_0\|_{L^2}+\|\sqrt{\rho_0}\theta_0\|_{L^2}
+\|\nabla u_0\|_{H^1}+\|\nabla \theta_0\|_{H^1}\notag\\
&+\|H_0\|_{H^2}
+ \|H_0\bar{x}^{\frac{a}{2}}\|_{H^1}+\|\rho_0\bar{x}^{a}\|_{L^{1}\cap H^{1}\cap W^{1,q}}+\|g_1\|_{L^2}+\|g_2\|_{L^2}.
\end{align*}
\end{proposition}

To show Proposition \ref{pro}, whose proof will be postponed to the end of this section, we begin with the following elementary estimate of the solutions.
\begin{lemma} \label{l3.01}
Let $(\rho,u,\theta,H)$ be as in Proposition \ref{pro} and $a>1$ be as in Theorem \ref{t1}. Then there exists a $T_1=T_1(N_0,E_0)>0$ such that, for all $t\in (0, T_1]$,
\begin{align}\label{igj1}
\sup_{0\le s\le t} \big(\|\sqrt{\rho}u\|_{L^2}^2+\|H\|_{L^2}^2+\|\rho\bar x^a\|_{L^1}+\|H\bar x^{\frac{a}{2}}\|_{L^2}^2\big)+\int_{0}^{t}\|\nabla H \bar x^{\frac{a}{2}}\|_{L^2}^2 ds\le C\exp\left\{C\int_{0}^{t}\psi^\alpha ds\right\},
\end{align}
where (and in what follows) $C$ denotes a generic positive constant depending on $\mu, \lambda, \nu, c_v, \kappa, q$, $a$, $\eta_0$, $N_0,$ and $E_0$, but independent of $R$, and we use $\alpha>1$ to denote a generic constant, which may be different from line to line.
\end{lemma}
\textit{Proof.}
1. First of all, the basic energy estimate gives that
\begin{equation}\label{3.1.0}
\int\Big(\frac12\rho|u|^2+c_v\rho\theta+\frac12|H|^2\Big)dx
=\int\Big(\frac12\rho_0|u_0|^2+c_v\rho_0\theta_0+\frac12|H_0|^2\Big)dx,
\quad t\geq0.
\end{equation}
Next, for $N>1$, let $\varphi_N\in C^\infty_0(B_N)$ satisfy
\begin{align}\label{vp1}
0\le \varphi_N \le 1, \ \varphi_N(x)=1,\ \mbox{if}\ |x|\le \frac{N}{2},\ \text{and} \ |\nabla \varphi_N|\le 3N^{-1}.
\end{align}
It follows from \eqref{mhd}$_1$ and  \eqref{3.1.0} that
\begin{align}\label{oo0}
\frac{d}{dt}\int \rho\varphi_{2N_0}dx
& =\int\rho u\cdot\nabla\varphi_{2N_0}dx
\ge -CN_0^{-1}\|\rho\|_{L^1}^{\frac12}
\|\sqrt{\rho}u\|_{L^2}\ge -\widetilde{C}(E_0),
\end{align}
where in the last inequality we have used
\begin{align*}
\int \rho dx=\int\rho_0dx,
\end{align*}
due to \eqref{mhd}$_1$. Integrating \eqref{oo0} and using \eqref{w1} give rise to
\begin{align}\label{p1}
\inf\limits_{0\le t\le T_1}\int_{B_{2N_0}} \rho dx
\ge \inf\limits_{0\le t\le T_1}\int \rho\varphi_{2N_0} dx
\ge \int \rho_0\varphi_{2N_0} dx-\widetilde{C}T_1
\ge \frac14,
\end{align}
where $T_1\triangleq\min\{1,(4\widetilde{C})^{-1}\}$. From now on, we will always assume that $t\le T_1.$
The combination of \eqref{p1}, \eqref{3.1.0}, and \eqref{3.i2} implies that, for $\varepsilon,\eta>0$ and $v\in \widetilde{D}^{1,2}(B_R)$ with $\sqrt{\rho}v\in L^2(B_R)$,
\begin{align}\label{3.v2}
\|v\bar x^{-\eta}\|_{L^{\frac{2+\varepsilon}{\tilde{\eta}}}}^2
&\le C(\varepsilon,\eta)\|\sqrt{\rho}v\|_{L^2}^2
+C(\varepsilon,\eta)\big(1+\|\rho\|_{L^\infty}^2\big)\|\nabla v\|_{L^2}^2,
\end{align}
where $\tilde{\eta}\triangleq\min\{1,\eta\}.$
It follows from H{\"o}lder's inequality, \eqref{3.2}, and \eqref{3.v2} that, for any $\varepsilon,\eta>0$ and $v\in \widetilde{D}^{1,2}(B_R)$ with $\sqrt{\rho}v\in L^2(B_R)$,
\begin{align}\label{local1}
\|\rho^\eta v\|_{L^{\frac{2+\varepsilon}{\tilde{\eta}}}}
& \le C\|\rho^\eta \bar x^{\frac{3\tilde{\eta} a}{4(2+\varepsilon)}} \|_{L^{ \frac{4(2+\varepsilon)}{3\tilde{\eta}}}} \|v\bar x^{-\frac{3\tilde{\eta} a}{4(2+\varepsilon)}} \|_{L^{ \frac{4(2+\varepsilon)}{\tilde{\eta}}}} \notag \\
& \le C\left(\int \rho^{\frac{4(2+\varepsilon)\eta}{3\tilde{\eta}}-1}\rho\bar x^a dx\right)^{ \frac{3\tilde{\eta}}{4(2+\varepsilon)}} \|v\bar x^{-\frac{3\tilde{\eta} a}{4(2+\varepsilon)}} \|_{L^{ \frac{4(2+\varepsilon)}{\tilde{\eta}}}} \notag \\
& \le C\|\rho\|_{L^\infty}^{\frac{4(2+\varepsilon)\eta-3\tilde{\eta}}{4(2+\varepsilon)}}
\|\rho\bar x^a\|_{L^1}^{\frac{3\tilde{\eta}}{4(2+\varepsilon)}}
\big(\|\sqrt{\rho}v\|_{L^2}+\big(1+\|\rho\|_{L^\infty}\big)\|\nabla v\|_{L^2}\big) \notag \\
& \le C\psi^\alpha\big(\|\sqrt{\rho}v\|_{L^2}+\psi^\alpha\|\nabla v\|_{L^2}\big),
\end{align}
where $\tilde{\eta}=\min\{1,\eta\}$. This together with \eqref{3.v2} implies that
\begin{align}\label{3.a2}
&\|\rho^\eta u\|_{L^{\frac{2+\varepsilon}{\tilde{\eta}}}}
+\|u\bar x^{-\eta}\|_{L^{\frac{2+\varepsilon}{\tilde{\eta}}}}
\le C\psi^\alpha\big(\|\sqrt{\rho}u\|_{L^2}+\psi^\alpha\|\nabla u\|_{L^2}\big)\le C\psi^\alpha, \\ \label{zx}
& \|\rho^\eta \theta\|_{L^{\frac{2+\varepsilon}{\tilde{\eta}}}}
+\|\theta\bar x^{-\eta}\|_{L^{\frac{2+\varepsilon}{\tilde{\eta}}}}
\le C\psi^\alpha\big(\|\sqrt{\rho}\theta\|_{L^2}+\psi^\alpha\|\nabla\theta\|_{L^2}\big)
\le C\psi^\alpha.
\end{align}

2. Noting that for any $\delta>0$, it holds that
\begin{align}\label{z10.6}
|\nabla\bar{x}|\leq C(\eta_0)\ln^{1+\eta_0}(e+|x|^2)
\leq C(\eta_0)\bar{x}^\delta.
\end{align}
Multiplying $\eqref{mhd}_{1}$ by $\bar x^a$ and integration by parts, we then obtain from H{\"o}lder's inequality, \eqref{2.5}, \eqref{3.v2}, \eqref{z10.6}, \eqref{3.a2}, and \eqref{3.1.0} that
\begin{align*}
\frac{d}{dt}\|\rho\bar x^a\|_{L^1}
& = \int\rho(u\cdot\nabla)\bar xa\bar x^{a-1}dx \notag \\
& \le C \int \rho|u|\bar x^{a-1+\frac{4}{8+a}}dx \notag \\
& \le C\|\rho\bar x^{a-1+\frac{8}{8+a}}\|_{L^\frac{8+a}{7+a}}
\|u\bar x^{-\frac{4}{8+a}}\|_{L^{8+a}} \notag \\
& \le C\psi^\alpha\|\rho\|_{L^\infty}^{\frac{1}{8+a}}\|\rho\bar x^a\|_{L^1}^{\frac{7+a}{8+a}}
\notag \\
& \le C\big(1+\|\rho\bar x^a\|_{L^1}\big)\psi^\alpha.
\end{align*}
This combined with Gronwall's inequality leads to
\begin{align}\label{igj1-2}
\sup_{0\le s\le t}\|\rho\bar x^a\|_{L^1}
\le C\exp\left\{C\int_{0}^{t}\psi^\alpha ds\right\}.
\end{align}

3. Multiplying \eqref{mhd}$_4$ by $H\bar{x}^a$ and integration by parts yield that
\begin{align}\label{lv4.1}
\frac{1}{2}\frac{d}{dt}\|H\bar{x}^{\frac{a}{2}}\|_{L^2}^2+\nu \|\nabla H \bar{x}^{\frac{a}{2}}\|_{L^2}^2
=& \frac{\nu}{2}\int |H|^2\Delta\bar{x}^adx+\int (H\cdot\nabla)u\cdot H\bar{x}^adx \notag\\
&-\frac{1}{2}\int \divv u |H|^2 \bar{x}^adx+\frac12\int |H|^2u\cdot\nabla\bar{x}^adx \triangleq \sum_{i=1}^4 \bar{I}_i.
\end{align}
Direct calculations lead to
\begin{equation}\label{10.4}
|\bar{I}_1|\leq C\int|H|^2\bar{x}^a\bar{x}^{-1}dx
\leq C\|H\bar{x}^{\frac{a}{2}}\|_{L^2}^2,
\end{equation}
and
\begin{align}\label{10.5}
|\bar{I}_2|+|\bar{I}_3| & \leq \int|\nabla u||H|^2\bar{x}^adx \notag \\
&\leq \|\nabla u\|_{L^2}\|H\bar{x}^{\frac{a}{2}}\|_{L^4}^2
 \notag \\
&\leq C\|\nabla u\|_{L^2}\|H\bar{x}^{\frac{a}{2}}\|_{L^2}\|H\bar{x}^{\frac{a}{2}}\|_{H^1}
 \notag \\
&\leq C\|\nabla u\|_{L^2}\|H\bar{x}^{\frac{a}{2}}\|_{L^2}
\big(\|H\bar{x}^{\frac{a}{2}}\|_{L^2}+\|\nabla H\bar{x}^{\frac{a}{2}}\|_{L^2}
+\|H\nabla\bar{x}^{\frac{a}{2}}\|_{L^2}\big) \notag \\
&\leq C\|\nabla u\|_{L^2}\|H\bar{x}^{\frac{a}{2}}\|_{L^2}
\big(\|H\bar{x}^{\frac{a}{2}}\|_{L^2}+\|\nabla H\bar{x}^{\frac{a}{2}}\|_{L^2}
+\|H\bar{x}^{\frac{a}{2}}\|_{L^2}
\|\bar{x}^{-1}\nabla\bar{x}\|_{L^\infty}\big)
 \notag \\
&\leq C\big(1+\|\nabla u\|_{L^2}^2\big)\|H\bar{x}^{\frac{a}{2}}\|_{L^2}^2
+\frac{\nu}{4}\|\nabla H\bar{x}^{\frac{a}{2}}\|_{L^2}^2\notag\\
&\le C\psi^\alpha\|H\bar{x}^{\frac{a}{2}}\|_{L^2}^2
+\frac{\nu}{4}\|\nabla H\bar{x}^{\frac{a}{2}}\|_{L^2}^2,
\end{align}
due to \eqref{2.4} and \eqref{z10.6}. Moreover,
it follows from H{\"o}lder's inequality, \eqref{z10.6}, \eqref{2.4}, and \eqref{3.v2} that
\begin{align}\label{10.6}
|\bar{I}_4| & \leq C\int|H|^2\bar{x}^a
\bar{x}^{-\frac12}|u|\bar{x}^{-\frac12+\frac15}dx \notag \\
&\leq C\|H\bar{x}^{\frac{a}{2}}\|_{L^4}
\|H\bar{x}^{\frac{a}{2}}\|_{L^2}
\|u\bar{x}^{-\frac34}\|_{L^4}\|\bar{x}^{-\frac{1}{20}}\|_{L^\infty}
\notag \\
&\leq C\Big(\|\sqrt{\rho}u\|_{L^2}^2+\big(1+\|\rho\|_{L^\infty}^2\big)\|\nabla u\|_{L^2}^2\Big)\|H\bar{x}^{\frac{a}{2}}\|_{L^2}^2
+\frac{\nu}{4}\|\nabla H\bar{x}^{\frac{a}{2}}\|_{L^2}^2\notag\\
&\leq C\psi^\alpha\|H\bar{x}^{\frac{a}{2}}\|_{L^2}^2
+\frac{\nu}{4}\|\nabla H\bar{x}^{\frac{a}{2}}\|_{L^2}^2.
\end{align}
Putting \eqref{10.4}--\eqref{10.6} into \eqref{lv4.1}, we thus deduce from Gronwall's inequality that
\begin{align}\label{lbqnew-gj10}
\sup_{0\le s\le t} \|H\bar{x}^{\frac{a}{2}}\|_{L^2}^2
+\int_{0}^{t}\|\nabla H\bar{x}^{\frac{a}{2}}\|_{L^2}^2ds\le C\exp\left\{C\int_{0}^{t}\psi^\alpha ds\right\}.
\end{align}
This along with \eqref{3.1.0} and \eqref{igj1-2} gives \eqref{igj1}.
\hfill $\Box$

\begin{lemma}\label{l3.0}
Let $(\rho,u,\theta,H)$ be as in Proposition \ref{pro} and $T_1$ be as in Lemma \ref{l3.01}. Then there exists a positive constant $\alpha>1$ such that, for all $t\in(0, T_1]$,
\begin{align}\label{gj3}
&\sup_{0\le s\le t}\big(\|\nabla u\|_{L^2}^2+\|H\|_{L^4}^4+\|\nabla H\|_{L^2}^2\big)
+\int_{0}^{t}\big(\|\sqrt{\rho}u_s\|_{L^2}^2+\| H_s\|_{L^2}^2+\|\Delta H\|_{L^2}^2\big)ds \notag \\
& \le C+C\int_{0}^{t}\psi^\alpha(s)ds.
\end{align}
\end{lemma}
\textit{Proof.}
1. Multiplying $\eqref{mhd}_{2}$ by $u_t$ and integration by parts show that
\begin{align}\label{3r1}
&\frac12\frac{d}{dt}\big(\mu\|\nabla u\|_{L^2}^2
+(\mu+\lambda)\|\divv u\|_{L^2}^2\big)
+\|\sqrt{\rho}u_t\|_{L^2}^2 \notag\\
&= -\int\rho u\cdot \nabla u\cdot u_tdx-\int\nabla P \cdot u_tdx
+\int \Big(H\cdot\nabla H-\frac{1}{2}\nabla|H|^2\Big)\cdot u_tdx\triangleq\sum_{i=1}^3 R_i.
\end{align}
By Cauchy-Schwarz inequality, H{\"o}lder's inequality, \eqref{2.4}, and \eqref{3.a2}, we get that
\begin{align}\label{cc5}
R_1=\left|-\int\rho u\cdot \nabla u\cdot u_tdx\right|
& \leq\frac12\int\rho|u_t|^2dx+\frac12\int\rho |u|^2|\nabla u|^2dx  \notag \\
& \le \frac12\|\sqrt{\rho}u_t\|_{L^2}^2
+\frac12\|\sqrt{\rho}u\|_{L^8}^{2}\|\nabla u\|_{L^{\frac83}}^{2} \notag \\
& \le \frac12\|\sqrt{\rho}u_t\|_{L^2}^2+C\|\sqrt{\rho}u\|_{L^8}^{2}\|\nabla u\|_{L^2}^{\frac32} \|\nabla u\|_{H^1}^{\frac12} \notag \\
& \le \frac12\|\sqrt{\rho}u_t\|_{L^2}^2+C\psi^{\alpha}.
\end{align}
For simplicity, we write $Q(\nabla u)=\frac{\mu}{2}|\nabla u+(\nabla u)^{tr}|^{2}+\lambda(\divv u)^2$, then integration by parts together with \eqref{mhd}$_3$, \eqref{3.a2}, \eqref{zx}, and Gagliardo-Nirenberg inequality indicates that
\begin{align}\label{gjr2}
R_2=&\int P\divv u_tdx\notag\\
=& \frac{d}{dt}\int P\divv udx-\int P_t\divv udx\notag\\
=& \frac{d}{dt}\int P\divv udx-\frac{R}{c_v}\int\big(Q(\nabla u)+\nu(\curl H)^2
+\kappa\Delta\theta-P\divv u-\frac{c_v}{R}\divv(P u)\big)\divv udx\notag\\
\le& \frac{d}{dt}\int P\divv udx+C\|\nabla u\|_{L^3}^3+C\|\nabla H\|_{L^3}^3+C\|\nabla^2\theta\|_{L^2}^2+C\|\nabla u\|_{L^2}^2+C\|\rho\|_{L^\infty}^{\frac{1}{2}}\|\sqrt{\rho}\theta\|_{L^2}\|\nabla u\|_{L^4}^2\notag\\
  &+C\|\rho\bar{x}^a\|_{W^{1,q}}\|\theta\bar{x}^{-\frac{a}{2}}\|_{L^{\frac{4q}{q-2}}}\|u\bar{x}^{-\frac{a}{2}}\|_{L^{\frac{4q}{q-2}}}\|\nabla u\|_{L^2}+C\|u\bar{x}^{-\frac{a}{2}}\|_{L^\infty}\|\rho\bar{x}^a\|_{L^\infty}\|\nabla \theta\|_{L^2}\|\nabla u\|_{L^2}\notag\\
\le& \frac{d}{dt}\int P\divv udx+C\psi^\alpha,
\end{align}
where in the last inequality one has used the following estimate
\begin{align*}
\|u\bar x^{-\frac{a}{2}}\|_{L^\infty}\le C\psi^\alpha.
\end{align*}
Indeed, it follows from \eqref{2.5}, Young's inequality, \eqref{z10.6}, \eqref{3.v2}, and \eqref{2.4} that
\begin{align}\label{3.22}
\|u\bar x^{-\frac{a}{2}}\|_{L^\infty}
& \le C\Big(\|u\bar x^{-\frac{a}{2}}\|_{L^4}
+\|\nabla(u\bar x^{-\frac{a}{2}})\|_{L^3}^{\frac35}
\|u\bar x^{-\frac{a}{2}}\|_{L^4}^{\frac25}\Big) \notag \\
& \le C\big(\|u\bar x^{-\frac{a}{2}}\|_{L^4}
+\|\nabla(u\bar x^{-\frac{a}{2}})\|_{L^3}\big) \notag \\
& \le C\big(\|u\bar x^{-\frac{a}{2}}\|_{L^4}
+\|u\bar x^{-\frac{a}{2}-1+\frac{a}{2}}\|_{L^3}
+\|\nabla u\|_{L^3}\big)
\notag \\
& \le C\Big(\|\sqrt{\rho}u\|_{L^2}+\big(1+\|\rho\|_{L^\infty}\big)\|\nabla u\|_{L^2}
+\|\nabla u\|_{L^2}^{\frac23}\|\nabla u\|_{H^1}^{\frac13}\Big)
\notag \\
& \le C\psi^\alpha.
\end{align}
Similarly, one has that
\begin{align}\label{3.3.1}
\|\theta\bar x^{-\frac{a}{2}}\|_{L^\infty}\le C\Big(\|\sqrt{\rho}\theta\|_{L^2}+\big(1+\|\rho\|_{L^\infty}\big)\|\nabla \theta\|_{L^2}
+\|\nabla\theta\|_{L^2}^{\frac23}\|\nabla\theta\|_{H^1}^{\frac13}\Big)\le C\psi^\alpha.
\end{align}
Integration by parts together with $\divv H=0$ and Gagliardo-Nirenberg inequality yields that
\begin{align}\label{gjr3}
R_3=&\int\Big(H\cdot\nabla H-\frac{1}{2}\nabla |H|^2\Big)\cdot u_tdx\notag\\
=&-\int H\cdot\nabla u_t\cdot Hdx+\frac{1}{2}\int|H|^2\divv u_tdx\notag\\
=&\frac{d}{dt}\Big(\frac{1}{2}\int|H|^2\divv udx-\int H\cdot\nabla u\cdot Hdx\Big)-\int H\cdot H_t \divv udx\notag\\
&+\int H_t\cdot\nabla u\cdot Hdx+\int H\cdot\nabla u\cdot H_tdx\notag\\
\le&\frac{d}{dt}\Big(\frac{1}{2}\int|H|^2\divv udx-\int H\cdot\nabla u\cdot Hdx\Big)+C\int |H|| H_t ||\nabla u|dx\notag\\
\le&\frac{d}{dt}\Big(\frac{1}{2}\int|H|^2\divv udx-\int H\cdot\nabla u\cdot Hdx\Big)+\delta \|H_t\|_{L^2}^2+C(\delta)\|H\|_{L^4}^2\|\nabla u\|_{L^4}^2\notag\\
\le&\frac{d}{dt}\Big(\frac{1}{2}\int|H|^2\divv udx-\int H\cdot\nabla u\cdot Hdx\Big)+\delta \|H_t\|_{L^2}^2+C(\delta)\|H\|_{L^2}\|H\|_{H^1}\|\nabla u\|_{L^2}\|\nabla u\|_{H^1}\notag\\
\le&\frac{d}{dt}\Big(\frac{1}{2}\int|H|^2\divv udx-\int H\cdot\nabla u\cdot Hdx\Big)+\delta \|H_t\|_{L^2}^2+C\psi^\alpha.
\end{align}
Inserting \eqref{cc5}, \eqref{gjr2}, and \eqref{gjr3} into \eqref{3r1} gives rise to
\begin{align}\label{3.8-1}
\frac{\mu}{2}\frac{d}{dt}\|\nabla u\|_{L^2}^2+\frac{1}{2}\|\sqrt{\rho}u_t\|_{L^2}^2\le B'(t)+\delta\|H_t\|_{L^2}^2+C\psi^\alpha,
\end{align}
where
$$B(t)\triangleq\int P\divv udx+\frac{1}{2}\int|H|^2\divv udx-\int H\cdot\nabla u\cdot Hdx$$
satisfies
\begin{align}\label{lv4.8'}
B(t)\le C\|H\|_{L^4}^4+C\|P\|_{L^2}^2+\frac{1}{2}\|\nabla u\|_{L^2}^2.
\end{align}
We deduce from \eqref{mhd}$_4$, \eqref{3.22}, and Gagliardo-Nirenberg inequality that
\begin{align}\label{gjht}
\|H_t\|_{L^2}^2& \leq C\big(\|\Delta H\|_{L^2}^2+\|u\cdot \nabla H\|_{L^2}^2+\|H\cdot \nabla u\|_{L^2}^2+\|H\divv u\|_{L^2}^2\big) \notag \\
& \le C\big(\|\nabla^2 H\|_{L^2}^2+\|u \bar x^{-\frac{a}{2}}\|_{L^{\infty}}^2\|\nabla H \bar x^{\frac{a}{2}}\|_{L^2}^2
+\|H\|_{L^4}^2\|\nabla u \|_{L^4}^2\big) \notag \\
& \le C\big(\|\nabla^2 H\|_{L^2}^2+\|u \bar x^{-\frac{a}{2}}\|_{L^{\infty}}^2\|\nabla H \bar x^{\frac{a}{2}}\|_{L^2}^2
+\|H \|_{L^2}\|H\|_{H^1}\|\nabla u \|_{L^2}\|\nabla u \|_{H^1}\big)\notag\\
&\le C\psi^\alpha,
\end{align}
which implies that
\begin{align}\label{xz3}
\int_{0}^{t}\|H_s\|_{L^2}^2ds \le C\int_{0}^{t}\psi^\alpha(s)ds.
\end{align}
Thus, integrating \eqref{3.8-1} over $[0,t]$ together with \eqref{lv4.8'} and \eqref{xz3} leads to
\begin{align}\label{gj3-1}
\sup_{0\le s\le t}\|\nabla u\|_{L^2}^2+\int_{0}^{t}\|\sqrt{\rho}u_s\|_{L^2}^2ds\le C+C\int_{0}^{t}\psi^\alpha(s)ds+C\|H\|_{L^4}^4+C\|P\|_{L^2}^2.
\end{align}

2. Multiplying\eqref{mhd}$_4$ by $|H|^2 H$ and integration by parts, we derive from \eqref{2.4} and \eqref{3.1.0} that
\begin{align}\label{3.1.7}
\frac{1}{4}\frac{d}{dt}\|H\|_{L^4}^4
+\nu\||\nabla H||H|\|_{L^2}^2+\frac{\nu}{2}\|\nabla|H|^2\|_{L^2}^2
& \leq C\|\nabla u\|_{L^2}\||H|^2\|_{L^4}^2 \notag \\
&\leq C\|\nabla u\|_{L^2}\||H|^2\|_{L^2}\||H|^2\|_{H^1} \notag \\
&\leq \frac{\nu}{4}\|\nabla |H|^2\|_{L^2}^2+C\|\nabla u\|_{L^2}^2\|H\|_{L^2}^2\|\nabla H\|_{L^2}^2\notag\\
&\le\frac{\nu}{4}\|\nabla |H|^2\|_{L^2}^2+C\psi^\alpha.
\end{align}
Integrating \eqref{3.1.7} over $[0,t]$ yields that
\begin{align}\label{3.1.6}
\sup_{0\le s\le t}\|H\|_{L^4}^4+\int_{0}^{t}\|\|\nabla H||H|\|_{L^2}^2ds \le C+C\int_{0}^{t}\psi^\alpha(s)ds.
\end{align}

3. It follows from \eqref{2.5} and \eqref{3.22} that
\begin{align}\label{3.53}
& \|\sqrt{\rho} u\|_{L^\infty} \notag \\
& \le\|\rho \bar{x}^a\|_{L^{\infty}}^\frac12\|u\bar{x}^{-\frac{a}{2}}\|_{L^\infty} \notag \\
& \le C\Big(\|\rho \bar{x}^a\|_{L^2}
+\|\rho \bar{x}^a\|^{\frac{q-2}{2(q-1)}}_{L^2}
\|\nabla(\rho \bar{x}^a)\|^{\frac{q}{2(q-1)}}_{L^q}\Big)^\frac12
\big(\|\sqrt{\rho}u\|_{L^2}+\big(1+\|\rho\|_{L^\infty}\big)\|\nabla u\|_{L^2}+\|\nabla u\|_{H^1}\big)\notag\\
&\le C\psi^\alpha,
\end{align}
Similarly, we have that
\begin{align}\label{3.3.2}
\|\sqrt{\rho} \theta\|_{L^\infty}\le C\psi^\alpha.
\end{align}
Multiplying \eqref{mhd}$_3$ by $P$ and integration by parts, we infer from H{\"o}lder's inequality and \eqref{3.3.2} that
\begin{align*}
\frac{1}{2}\frac{d}{dt}\|P\|_{L^2}^2&\le C\int P^2 |\divv u|dx+C\int P(|\nabla u|^2+|\nabla H|^2)dx+C\int P|\Delta\theta| dx\notag\\
&\le C\|\rho\|_{L^\infty}\|\sqrt{\rho}\theta\|_{L^\infty}\|\sqrt{\rho}\theta\|_{L^2}\|\nabla u\|_{L^2}+ C\|\rho\|_{L^\infty}^{\frac{1}{2}}\|\sqrt{\rho}\theta\|_{L^\infty}(\|\nabla u\|_{L^2}^2+\|\nabla H\|_{L^2}^2)\notag\\
&\quad+C\|\rho\|_{L^\infty}^\frac{1}{2}\|\sqrt{\rho} \theta\|_{L^2}\|\nabla^2\theta\|_{L^2}
\le C\psi ^\alpha,
\end{align*}
which yields that
\begin{align}\label{3.1.8}
\sup_{0\le s\le t}\|P\|_{L^2}^2\le C+C\int_0^t \psi^\alpha ds.
\end{align}

4. Multiplying \eqref{mhd}$_4$ by $\Delta H$ and integration by parts, we get from Gagliardo-Nirenberg inequality that
\begin{align*}
&\frac{d}{dt}\int |\nabla H|^2 dx +2\nu \int|\Delta H|^2dx \notag \\
&\le C\int|\nabla u||\nabla H|^2dx+ C\int|\nabla u||H||\Delta H|dx\notag\\
&\le C\|\nabla u\|_{L^3}\|\nabla H\|_{L^2}^\frac{4}{3}\|\nabla H\|_{H^1}^\frac{2}{3}+C\|\nabla u\|_{L^3}\|H\|_{L^6}\|\Delta H\|_{L^2}\notag\\
&\le C\|\nabla u\|_{L^2}^2\|\nabla u\|_{H^1}+C\|\nabla H\|_{L^2}^4+C\|\nabla ^2 H\|_{L^2}^2+C\|H\|_{L^2}^2\|H\|_{H^1}^4
\le C\psi^\alpha,
\end{align*}
which implies that
 \begin{align}\label{3.1.10}
\sup_{0\le s\le t}\|\nabla H\|_{L^2}^2+\int_0^t\|\Delta H\|_{L^2}^2ds
\le C+C\int_{0}^{t}\psi^\alpha(s)ds.
\end{align}
Hence, the desired \eqref{gj3} follows directly from \eqref{xz3}, \eqref{gj3-1}, \eqref{3.1.6}, \eqref{3.1.8}, and \eqref{3.1.10}.
\hfill $\Box$

\begin{lemma}\label{le-3}
Let $(\rho,u,\theta,H)$ be as in Proposition \ref{pro} and $T_1$ be as in Lemma \ref{l3.01}. Then there exists a positive constant $\alpha>1$ such that, for all $t\in (0, T_1],$
\begin{align}\label{gj6}
\sup_{0\le s\le t}\big(\|\sqrt{\rho}u_s\|_{L^2}^2+\||H||\nabla H|\|_{L^2}^2\big)+\int_0^t(\|\nabla u_s\|_{L^2}^2+\||\Delta H||H|\|_{L^2}^2)ds \le C\exp\left\{C\int_0^t\psi^\alpha ds\right\}.
\end{align}
\end{lemma}
\textit{Proof.}
1. Differentiating $\eqref{mhd}_2$ with respect to $t$ gives that
\begin{align}\label{zb1}
&\rho u_{tt}+\rho u\cdot \nabla u_t-\mu\Delta u_t-(\mu+\lambda)\nabla\divv u_t\notag\\
&=-\rho_t(u_t+u\cdot\nabla u)-\rho u_t\cdot\nabla u -\nabla P_t
+\Big(H\cdot\nabla H-\frac12\nabla |H|^2\Big)_t.
\end{align}
Multiplying \eqref{zb1} by $u_t$ and integrating the resulting equality by parts over $B_R$, we obtain after using $\eqref{mhd}_1$ and $\divv H=0$ that
\begin{align}\label{na8}
&\frac{1}{2}\frac{d}{dt} \|\sqrt{\rho}u_t\|_{L^2}^2+\mu\|\nabla u_t\|_{L^2}^2 +(\mu+\lambda)\|\divv u_t\|_{L^2}^2 \notag\\
&\le-\int\rho_t|u_t|^2dx-\int(\rho u)_t\cdot\nabla u\cdot u_tdx-\int u_t\cdot\nabla P_tdx\notag\\
&\quad+C\int H_t\cdot \nabla H\cdot u_tdx+C\int H\cdot \nabla H_t\cdot u_t dx\triangleq \sum_{i=1}^5 \hat{I}_i.
\end{align}
Using \eqref{mhd}$_1$, \eqref{3.i2}, \eqref{3.53}, \eqref{3.3.2}, and Gagliardo-Nirenberg inequality, one deduces that
\begin{align}\label{na2}
\hat{I}_1\le & C\int \rho|u||u_t||\nabla u_t|dx
\le  \frac{\mu}{8}\|\nabla u_{t}\|_{L^2}^{2}+C\|\sqrt{\rho}u\|_{L^\infty}^2\|\sqrt{\rho}u_t\|_{L^2}^2 \le  \frac{\mu}{8}\|\nabla u_{t}\|_{L^2}^{2}+C\psi^{\alpha},\\ \label{5.a3}
\hat{I}_2\le&C\int(|u||\nabla \rho|+\rho \divv u)|u||\nabla u ||u_t|dx+C\int\rho |u_t||\nabla u|| u_t|dx\notag\\
\le&C\|\rho \bar{x}^a\|_{W^{1,q}}\|\nabla u\|_{L^2}\|u\bar{x}^{-\frac{a}{3}}\|_{\frac{6q}{q-2}}^2\|u_t\bar{x}^{-\frac{a}{3}}\|_{\frac{6q}{q-2}}+C\|\sqrt{\rho}u\|_{L^\infty}\|\nabla u\|_{L^4}^2\|\sqrt{\rho}u_t\|_{L^2}\notag\\
&+C\|\rho \bar{x}^a\|_{L^\infty}^{\frac{1}{2}}\|\nabla u\|_{L^4}\|\sqrt{\rho}u_t\|_{L^2}\|u_t\bar{x}^{-\frac{a}{2}}\|_{L^4}\notag\\
\le&\frac{\mu}{8} \| \nabla u_{t}\|_{L^{2}}^{2} + C \psi^{\alpha},\\ \label{3.3.3}
\hat{I}_3\le&C\int\big(|u||\nabla \rho|+\rho \divv u\big)\theta \divv u_tdx+C\int\rho \theta_t\divv u_tdx\notag\\
\le& C\|\rho \bar{x}^a\|_{W^{1,q}}\|u\bar{x}^{-\frac{a}{2}}\|_{\frac{4q}{q-2}}\|\theta\bar{x}^{-\frac{a}{2}}\|_{\frac{4q}{q-2}}\|\nabla u_t\|_{L^2}+C\big(\|\rho \theta\|_{L^\infty}\|\nabla u\|_{L^2}+\|\rho\|_{L^\infty}^{\frac{1}{2}}\|\sqrt{\rho} \theta_t\|_{L^2}\big)\|\nabla u_t\|_{L^2}\notag\\
\le& \frac{\mu}{8} \| \nabla u_{t}\|_{L^{2}}^{2} + C \psi^{\alpha}.
\end{align}
It follows from integration by parts, H\"older's inequality, \eqref{gjht}, \eqref{mhd}$_4$, and \eqref{2.4} indicates that
\begin{align}\label{na3}
\hat{I}_4+\hat{I}_5
& =-C\int H_t\cdot\nabla u_t\cdot H dx-C\int H\cdot\nabla u_t\cdot H_tdx \notag \\
& \leq C\int|H_t||H||\nabla u_t|dx\le C\int(|\Delta H|+|\nabla u||H|+|\nabla H||u|)|H||\nabla u_t| dx\notag \\
&\le \frac{\mu}{8}\|\nabla u_{t}\|_{L^2}^2+C\big(\||\Delta H||H|\|_{L^2}^2+\||\nabla u||H|^2\|_{L^2}^2+\||\nabla H||u||H|\|_{L^2}^2\big)
\notag \\
&\le \frac{\mu}{8}\|\nabla u_{t}\|_{L^2}^2+C\||\Delta H||H|\|_{L^2}^2+C\||\nabla u\|_{L^4}^2\||H|^2\|_{L^4}^2\notag\\
&\quad +C\||u|^2\bar{x}^{-\frac{a}{2}}\|_{L^8}\|H\bar{x}^{\frac{a}{2}}\|_{L^2}\|H\|_{L^8}\||\nabla H|^2\|_{L^4}
\notag \\
& \le \frac{\mu}{8}\|\nabla u_{t}\|_{L^2}^2+C\||\Delta H||H|\|_{L^2}^2+C\||\nabla u\|_{L^2}\||\nabla u\|_{H^1}\||H|^2\|_{L^2}\|\nabla|H|^2\|_{L^2}\notag\\
&\quad +C\|u\bar{x}^{-\frac{a}{4}}\|_{L^{16}}^2\|H\bar{x}^{\frac{a}{2}}\|_{L^2}\|H\|_{L^2}^{\frac{1}{4}}\|\nabla H\|_{L^2}^{\frac{3}{4}}\|\nabla H\|_{L^2}^{\frac{1}{2}}\|\nabla H\|_{H^1}^{\frac{3}{2}}\notag \\
& \le \frac{\mu}{8}\|\nabla u_{t}\|_{L^2}^2+C\||\Delta H||H|\|_{L^2}^2+C\psi^\alpha\|\nabla |H|^2\|_{L^2}^2+C\psi^\alpha
\notag \\
& \le \frac{\mu}{8}\|\nabla u_{t}\|_{L^2}^2+C_1\||\Delta H||H|\|_{L^2}^2+C\psi^\alpha\big(1+\||\nabla H||H|\|_{L^2}^2\big).
\end{align}
Thus, substituting \eqref{na2}--\eqref{na3} into \eqref{na8}, we obtain that
\begin{align}\label{a4.6}
\frac{d}{dt}\|\sqrt{\rho}u_{t}\|_{L^2}^2+\mu\|\nabla u_{t}\|_{L^2}^2
\le C_1\||\Delta H||H|\|_{L^2}^2+C\psi^\alpha\big(1+\||\nabla H||H|\|_{L^2}^2\big).
\end{align}

2. We will use the methods of \cite{LB2015} to control the right-hand term of \eqref{a4.6}. For $a_1,a_2\in\{-1,0,1\}$, denoting
\begin{align}
\tilde{H}(a_1,a_2)=a_1H^1+a_2H^2,\,\,\,\,\tilde{u}(a_1,a_2)=a_1u^1+a_2u^2,
\end{align}
it thus follows from \eqref{mhd}$_4$ that
\begin{align}\label{3.2.7}
\tilde{H}_t-\nu\Delta\tilde{H}=H\cdot\nabla\tilde{u}-u\cdot\nabla\tilde{H}-\tilde{H}\divv u.
\end{align}
Multiplying \eqref{3.2.7} by $4\nu^{-1}\tilde{H}\Delta|\tilde{H}|^2$ and integrating the resulting equations by parts lead to
\begin{align}\label{3.2.8}
\nu^{-1}(\|\nabla |\tilde{H}|^2\|_{L^2}^2)_t+2\|\Delta|\tilde{H}|^2\|_{L^2}^2&=4\int|\nabla\tilde{H}|^2\Delta|\tilde{H}|^2dx-4\nu^{-1}\int H\cdot\nabla\tilde{u}\tilde{H}\Delta|\tilde{H}|^2dx\notag\\
& \quad +4\nu^{-1}\int\divv u|\tilde{H}|^2\Delta|\tilde{H}|^2dx+2\nu^{-1}\int u\cdot\nabla|\tilde{H}|^2\Delta|\tilde{H}|^2dx\notag\\
& \le C\|\nabla u\|_{L^4}^4+C\|\nabla H\|_{L^4}^4+C\||H|^2\|_{L^4}^4+\|\Delta |\tilde{H}|^2\|_{L^2}^2\notag\\
& \le C\psi^\alpha\big(1+\||\nabla H||H|\|_{L^2}^2\big)+\|\Delta |\tilde{H}|^2\|_{L^2}^2,
\end{align}
where in the first inequality we have used the following estimate
\begin{align*}
\int u\cdot\nabla|\tilde{H}|^2\Delta|\tilde{H}|^2&=-\int\nabla u\cdot\nabla|\tilde{H}|^2\cdot\nabla|\tilde{H}|^2dx+\frac{1}{2}\int\divv u|\nabla |\tilde{H}|^2|^2dx\notag\\
&\le C\|\nabla u\|_{L^4}^4+C\|\nabla H\|_{L^4}^4+C\||H|^2\|_{L^4}^4.
\end{align*}
Noticing that
\begin{align}\label{3.2.9}
\||\nabla H||H|\|_{L^2}^2
\le\|\nabla|\tilde{H}(1,0)|^2\|_{L^2}^2+\|\nabla|\tilde{H}(0,1)|^2\|_{L^2}^2
+\|\nabla|\tilde{H}(1,1)|^2\|_{L^2}^2+\|\nabla|\tilde{H}(1,-1)|^2\|_{L^2}^2,
\end{align}
and
\begin{align}
\||\Delta H||H|\|_{L^2}^2
&\le C\|\nabla H\|_{L^4}^4+C\|\Delta|\tilde{H}(1,0)|^2\|_{L^2}^2
+\|\Delta|\tilde{H}(0,1)|^2\|_{L^2}^2 \notag \\
& \quad +\|\Delta|\tilde{H}(1,1)|^2\|_{L^2}^2
+\|\Delta|\tilde{H}(1,-1)|^2\|_{L^2}^2,
\end{align}
adding \eqref{3.2.8} multiplied by $4(C_1+1)$ to \eqref{a4.6}, we deduce from  Gagliardo-Nirenberg inequality and \eqref{3.2.9} that
\begin{align}\label{3.2.12}
\frac{d}{dt}\big(\|\sqrt{\rho} u_t\|_{L^2}^2+\||\nabla H||H|\|_{L^2}^2\big)+\mu\|\nabla u_t\|_{L^2}^2+\||\Delta H||H|\|_{L^2}^2\le C\psi^\alpha\big(1+\||\nabla H||H|\|_{L^2}^2\big).
\end{align}

3. We derive from \eqref{mhd}$_2$, \eqref{3.53}, and \eqref{c.c} that
\begin{align}\label{xz4}
\int \rho|u_t|^2(x,0)dx
& \le\lim_{t \rightarrow 0}\sup\int\rho^{-1}|\mu\Delta u+(\mu+\lambda)\nabla\divv u+H\cdot\nabla H-\frac{1}{2}\nabla |H|^2
-\nabla P-\rho u\cdot\nabla u|^2dx \notag \\
& \le 2\|g_1\|_{L^2}^2+2\|\rho u(0)\|_{L^\infty}^2\|\nabla u(0)\|_{L^2}^2 \notag \\
& \leq 2\|g_1\|_{L^2}^2+CE_0^\alpha \le C.
\end{align}
Moreover, it follows from \eqref{2.4} and H{\"o}lder's inequality that
\begin{align*}
\||H_0||\nabla H_0|\|_{L^2}^2\le C\|H_0\|_{L^4}^2\|\nabla H_0\|_{L^4}^2\le C\|H_0\|_{L^2}\|H_0\|_{H^1}\|\nabla H_0\|_{L^2}\|\nabla H_0\|_{H^1}
\le CE_0^\alpha\le C.
\end{align*}
This along with \eqref{3.2.12}, Gronwall's inequality, and \eqref{xz4} leads to \eqref{gj6}.\hfill $\Box$

\begin{lemma} \label{le-3'}
Let $(\rho,u,\theta,H)$ be as in Proposition \ref{pro} and $T_1$ be as in Lemma \ref{l3.01}. Then there exists a positive constant $\alpha>1$ such that, for all $t\in(0, T_1],$
\begin{align}\label{gj10'}
&\sup_{0\le s\le t}\|\nabla H\bar{x}^{\frac{a}{2}}\|_{L^2}^2+\int_{0}^{t}\Big(\|\nabla^2u\|_{ L^q}^{\frac{q+1}{q}}
+\|\nabla^2u\|_{L^q}^2
+\|\nabla^2H\bar{x}^{\frac{a}{2}}\|_{L^2}^2\Big)ds \notag \\
& \le C\exp\left\{C\exp\left\{C\int_{0}^{t} \psi^{\alpha} ds\right\}\right\}.
\end{align}
\end{lemma}
\textit{Proof.}
1. Multiplying \eqref{mhd}$_4$ by $\Delta H\bar{x}^a$ and integration by parts lead to
\begin{align}\label{AMSS5}
&\frac{1}{2}\frac{d}{dt}\int |\nabla H|^2\bar{x}^adx+\nu \int |\Delta H|^2\bar{x}^adx \notag \\
& \le C\int|\nabla H| |H| |\nabla u| |\nabla\bar{x}^a|dx+C\int|\nabla H|^2|u| |\nabla\bar{x}^a|dx+C\int|\nabla H| |\Delta H| | \bar{x}^a|dx \notag \\
& \quad+ C\int |H||\nabla u||\Delta H|\bar{x}^adx+C\int |\nabla u||\nabla H|^2 \bar{x}^adx
\triangleq \sum_{i=1}^5 J_i.
\end{align}
Using \eqref{lbqnew-gj10}, \eqref{3.v2}, H\"older's inequality, and \eqref{2.4}, we get by some direct calculations that
\begin{align*}
J_1\le & C\|H\bar{x}^{\frac{a}{2}}\|_{L^4}\|\nabla u\|_{L^4}\|\nabla H\bar{x}^{\frac{a}{2}}\|_{L^2}\\
\le & C\|H\bar{x}^{\frac{a}{2}}\|_{L^2}^{\frac12}
\big(\|\nabla H\bar{x}^{\frac{a}{2}}\|_{L^2}
+\|H\bar{x}^{\frac{a}{2}}\|_{L^2}\big)^{\frac12}
\|\nabla u\|_{L^2}^{\frac12}\|\nabla u\|_{H^1}^{\frac12}
\|\nabla H\bar{x}^{\frac{a}{2}}\|_{L^2}\\
\le & C\psi^{\alpha}+C\psi^{\alpha}\|\nabla H\bar{x}^{\frac{a}{2}}\|_{L^2}^2,\\
J_2\leq & C\||\nabla H|^{2-\frac{2}{3a}}\bar{x}^{a-\frac13}\|_{L^{\frac{6a}{6a-2}}} \|u\bar{x}^{-\frac13}\|_{L^{6a}}\||\nabla H|^{\frac{2}{3a} }\|_{L^{6a}}\\
\le & C\psi^{\alpha}\|\nabla H\bar{x}^{\frac{a}{2}} \|_{L^2}^\frac{6a-2}{3a}\|\nabla H\|_{L^4}^\frac{2}{3a} \\
\leq & C\psi^{\alpha}\|\nabla H\bar{x}^{\frac{a}{2}}\|_{L^2}^2
+C\|\nabla H\|_{L^4}^2\\
\leq & C\psi^{\alpha}\|\nabla H \bar{x}^{\frac{a}{2}}\|_{L^2}^2
+\frac{\nu}{4}\|\Delta H \bar{x}^{\frac{a}{2}}\|_{L^2}^2,\\
J_3+J_4\le &\frac{\nu}{4}\|\Delta H\bar{x}^{\frac{a}{2}}\|_{L^2}^2
+C\|\nabla H\bar{x}^{\frac{a}{2}}\|_{L^2}^2
+C\|H\bar{x}^{\frac{a}{2}}\|_{L^4}^2 \|\nabla u\|_{L^4}^2\\
\le &\frac{\nu}{4}\|\Delta H\bar{x}^{\frac{a}{2}}\|_{L^2}^2+C\|\nabla H\bar{x}^{\frac{a}{2}}\|_{L^2}^2+C\|H\bar{x}^{\frac{a}{2}}\|_{L^2}\big(\|\nabla H\bar{x}^{\frac{a}{2}}\|_{L^2}+\|H\bar{x}^{\frac{a}{2}}\|_{L^2}\big)
\|\nabla u\|_{L^2} \|\nabla u\|_{H^1} \\
\le & \frac{\nu}{4}\|\Delta H\bar{x}^{\frac{a}{2}}\|_{L^2}^2+C\psi^{\alpha}\|\nabla H\bar{x}^{\frac{a}{2}}\|_{L^2}^2+C\psi^{\alpha},\\
J_5 \le & C\|\nabla u\|_{L^\infty} \|\nabla H\bar{x}^{\frac{a}{2}}\|_{L^2}^2
\le  C\Big(\psi^{\alpha}+\|\nabla^2 u\|_{L^q}^{\frac{q+1}{q}}\Big)\|\nabla H\bar{x}^{\frac{a}{2}}\|_{L^2}^2.
\end{align*}
Hence, substituting the above estimates into \eqref{AMSS5}
and noting
\begin{align*}
\int|\nabla^2H|^2\bar{x}^adx
& = \int|\Delta H|^2\bar{x}^adx
-\int\partial_i\partial_kH\cdot\partial_kH\partial_i\bar{x}^adx
+\int\partial_i\partial_iH\cdot\partial_kH\partial_k\bar{x}^adx \\
& \leq \int|\Delta H|^2\bar{x}^adx+\frac12\int|\nabla^2H|^2\bar{x}^adx
+C\int|\nabla H|^2\bar{x}^adx,
\end{align*}
we derive that
\begin{align}\label{AMSS10}
& \frac{d}{dt}\|\nabla H\bar{x}^{\frac{a}{2}}\|_{L^2}^2
+\nu\|\nabla^2H\bar{x}^{\frac{a}{2}}\|_{L^2}^2
\le C\Big(\psi^{\alpha}+\|\nabla^2 u\|_{L^q}^{\frac{q+1}{q}}\Big)
\|\nabla H\bar{x}^{\frac{a}{2}}\|_{L^2}^2+C\psi^{\alpha}.
\end{align}
We now claim that
\begin{align}\label{gj7}
\int_{0}^{t}\Big(\|\nabla^2u\|_{ L^q}^{\frac{q+1}{q}}
+\|\nabla^2u\|_{L^q}^2
\Big)ds
\le C\exp\left\{C\int_0^t\psi^\alpha ds\right\},
\end{align}
whose proof will be given at the end of this proof. Thus, by \eqref{AMSS10}, \eqref{lbqnew-gj10}, \eqref{gj3}, \eqref{gj7}, and Gronwall's inequality,
\begin{align}\label{igj10'}
\sup_{0\le s\le t}\|\nabla H\bar{x}^{\frac{a}{2}}\|_{L^2}^2
+\int_{0}^{t}\|\nabla^2H\bar{x}^{\frac{a}{2}}\|_{L^2}^2ds
\le C\exp\left\{C\exp\left\{C\int_{0}^{t} \psi^{\alpha} ds\right\}\right\}.
\end{align}

2. We next show \eqref{gj7}. It follows from \eqref{mhd}$_2$, \eqref{lemma2.4}, \eqref{2.4}, and \eqref{3.3.1} that
\begin{align}\label{lv3.60}
\|\nabla^2u\|_{L^q}&\leq C\big(\|\rho u_t\|_{L^q}+\|\rho u\cdot\nabla u\|_{L^q}+\|\nabla P\|_{L^q}+\||H||\nabla H|\|_{L^q} \big)\notag \\
&\leq C\|\rho\|_{L^\infty}^{\frac{1}{2}}\|\sqrt{\rho}u_t\|_{L^2}^{\frac{2(q-1)}{q^2-2}}\|\sqrt{\rho}u_t\|_{L^{q^2}}^{\frac{q^2-2q}{q^2-2}}
+C\|\rho u\|_{L^{2q}}\|\nabla u\|_{L^{2q}}\notag\\
&\quad
+C\big(\|\nabla \rho \bar{x}^a\|_{L^q}\|\theta\bar{x}^{-a}\|_{L^\infty}+\|\rho\|_{L^\infty}\|\nabla \theta\|_{L^q}\big)
+C\|H\|_{L^{2q}}\|\nabla H\|_{L^{2q}} \notag \\
& \le C\psi^\alpha\|\sqrt{\rho}u_t\|_{L^2}^{\frac{2(q-1)}{q^2-2}}
\Big(\|\sqrt{\rho}u_t\|_{L^2}+\big(1+\|\rho\|_{L^\infty}\big)\|\nabla u_t\|_{L^2}\Big)^{\frac{q^2-2q}{q^2-2}}\notag\\
 &\quad+C\psi^\alpha\Big(\|\sqrt{\rho}u\|_{L^2}+\big(1+\|\rho\|_{L^\infty}\big)\|\nabla u\|_{L^2}\Big)
\|\nabla u\|_{L^2}^{\frac1q}\|\nabla u\|_{H^1}^{\frac{q-1}{q}} \notag \\
& \quad +C\psi^\alpha+C\psi^\alpha\|\nabla\theta\|_{L^2}^{\frac{2}{q}}\|\nabla\theta\|_{H^1}^\frac{q-2}{q}+C\|H\|_{L^2}^{\frac1q}\|H\|_{H^1}^{\frac{q-1}{q}}
\|\nabla H\|_{L^2}^{\frac1q}\|\nabla H\|_{H^1}^{\frac{q-1}{q}} \notag \\
&\le C\psi^\alpha
\|\nabla u_t\|_{L^{2}}^{\frac{q^2-2q}{q^2-2}}+C\psi^\alpha,
\end{align}
which together with Young's inequality and \eqref{gj6} implies that
\begin{align} \label{olv3.61}
\int_0^t\|\nabla^2u\|_{L^q}^{\frac{q+1}{q}}ds
& \le C\int_0^t\psi^\alpha\left(\|\nabla u_s\|^2_{L^2}\right)^{\frac{(q-2)(q+1)}{2(q^2-2)}} ds
+C\int_0^t \psi^\alpha ds \notag \\
& \le C\int_0^t\left(\psi^\alpha+\|\nabla u_s\|_{L^2}^2\right)ds
\le C\exp\left\{C\int_0^t\psi^\alpha ds\right\},
\end{align}
and
\begin{align}\label{olv3.62}
\int_0^t\|\nabla^2u\|_{L^q}^2ds
& \le C\int_0^t \psi^\alpha
\big(\|\nabla u_s\|_{L^{2}}^2\big)^{\frac{q^2-2q}{q^2-2}}ds
+C\int_0^t\psi^\alpha ds \notag \\
&\le C\int_0^t\|\nabla u_s\|_{L^{2}}^2ds+C\int_0^t\psi^\alpha ds
\le C\exp\left\{C\int_0^t\psi^\alpha ds\right\},
\end{align}
where we have used $\frac{(q-2)(q+1)}{2(q^2-2)},\frac{q^2-2q}{q^2-2}\in(0,1)$ due to $q>2$.
One obtains \eqref{gj7} from \eqref{olv3.61} and \eqref{olv3.62}.
\hfill $\Box$

\begin{lemma}\label{le-4}
Let $(\rho,u,\theta,H)$ be as in Proposition \ref{pro} and $T_1$ be as in Lemma \ref{l3.01}. Then there exists a positive constant $\alpha>1$ such that, for all $t\in (0, T_1]$,
\begin{align}\label{gj8}
&\sup\limits_{0\le s\le t} \|\rho\bar x^a \|_{ L^1\cap H^1\cap W^{1,q}} \le C\exp\left\{C\exp\left\{C\int_{0}^{t} \psi^{\alpha} ds\right\}\right\}.
\end{align}
\end{lemma}
\textit{Proof.} We derive from \eqref{mhd}$_1$ that $\rho\bar{x}^a$ satisfies
\begin{equation}\label{A}
\partial_{t}(\rho\bar{x}^a)+ u \cdot\nabla(\rho\bar{x}^a)
-a\rho\bar{x}^{a}u\cdot\nabla(\ln\bar{x})+\rho \bar{x}^a\divv u=0.
\end{equation}
Operating $\nabla$ to \eqref{A} and then multiplying the resultant equation by $|\nabla(\rho\bar{x}^a)|^{r-2}\nabla(\rho\bar{x}^a)$ for $r\in[2,q]$,
we obtain after integration by parts that
\begin{align}\label{6.4}
\frac{d}{dt}\|\nabla(\rho\bar{x}^a)\|_{L^r} \leq
& C\left(1+\|\nabla u \|_{L^\infty}
+\|u\cdot\nabla(\ln\bar{x})\|_{L^\infty}\right)
\|\nabla(\rho\bar{x}^a)\|_{L^r} \notag \\
& +C\|\rho\bar{x}^a\|_{L^\infty}\big(\||\nabla u ||\nabla(\ln\bar{x})|\|_{L^r}+\||u||\nabla^{2}(\ln\bar{x})|\|_{L^r}+\|\nabla^2u\|_{L^r}\big).
\end{align}
By \eqref{2.5} and Young's inequality, we see that
\begin{align*}
\|\nabla u\|_{L^\infty}\leq C\|\nabla u\|_{L^2}
+C\|\nabla u\|_{L^2}^{\frac{q-2}{2q-2}}\|\nabla^2u\|_{L^q}^{\frac{q}{2q-2}}
\leq C\|\nabla^2u\|_{L^q}^2+C\psi^\alpha.
\end{align*}
Similarly to \eqref{3.22}, we obtain from \eqref{z10.6} that
\begin{align*}
\|u\cdot\nabla(\ln\bar{x})\|_{L^\infty}
=\|u\cdot\bar{x}^{-1}\nabla\bar{x}\|_{L^\infty}
\leq C\psi^\alpha.
\end{align*}
From \eqref{2.5}, we have
\begin{align*}
\|\rho\bar{x}^a\|_{L^\infty}
\leq C\|\rho\bar{x}^a\|_{L^2}+C\|\rho\bar{x}^a\|_{L^2}^{\frac{q-2}{2q-2}}
\|\nabla(\rho\bar{x}^a)\|_{L^q}^{\frac{q}{2q-2}}\leq C\psi^\alpha.
\end{align*}
Applying \eqref{z10.6} and \eqref{2.4}, we get that
\begin{align*}
\||\nabla u ||\nabla(\ln\bar{x})|\|_{L^r}
\leq C\|\nabla u\|_{L^r}\|\bar{x}^{-\frac{4+a}{8+a}}\|_{L^\infty}
\leq C\|\nabla u\|_{L^2}^{\frac2r}\|\nabla u\|_{H^1}^{\frac{r-2}{r}}
\leq C\psi^\alpha.
\end{align*}
Moreover, it follows from \eqref{z10.6} and \eqref{3.a2} that
\begin{align*}
\||u||\nabla^{2}(\ln\bar{x})|\|_{L^r}\leq C\psi^\alpha.
\end{align*}
As a consequence, inserting the above estimates into \eqref{6.4}, we derive that
\begin{align*}
\frac{d}{dt}\big(\|\nabla(\rho\bar{x}^a)\|_{L^2}+\|\nabla(\rho\bar{x}^a)\|_{L^q}\big)
\leq C\big(\psi^\alpha+\|\nabla^2 u\|_{L^q}^2\big)
\big(1+\|\nabla(\rho\bar{x}^a)\|_{L^2}+\|\nabla(\rho\bar{x}^a)\|_{L^q}\big),
\end{align*}
which combined with Gronwall's inequality, \eqref{gj10'}, and \eqref{olv3.62} indicates that
\begin{align}\label{3.68}
&\sup\limits_{0\le s\le t} \|\nabla(\rho\bar x^a)\|_{L^2\cap L^q}\le C\exp\left\{C\exp\left\{C\int_{0}^{t} \psi^{\alpha} ds\right\}\right\}.
\end{align}
Similarly, multiplying \eqref{A} by $(\rho\bar{x}^a)^{r-1}$ for $r\in[2,q]$ and then integrating the resultant equation over $B_R$,
we can also deduce that
\begin{align*}
&\sup\limits_{0\le s\le t} \|\rho\bar x^a\|_{L^2\cap L^q}\le C\exp\left\{C\exp\left\{C\int_{0}^{t} \psi^{\alpha} ds\right\}\right\}.
\end{align*}
This along with \eqref{igj1} and \eqref{3.68} implies \eqref{gj8}.
\hfill $\Box$

\begin{lemma}\label{le--3}
Let $(\rho,u,\theta,H)$ be as in Proposition \ref{pro} and $T_1$ be as in Lemma \ref{l3.01}. Then there exists a positive constant $\alpha>1$ such that, for all $t\in (0, T_1],$
\begin{align}\label{gj9'}
\sup_{0\le s\le t}\|H_s\|_{L^2}^2+\int_0^t\|\nabla H_s\|_{L^2}^2ds \le C\exp\left\{C\exp\left\{C\int_{0}^{t} \psi^{\alpha} ds\right\}\right\}.
\end{align}
\end{lemma}
\textit{Proof.}
1. Differentiating $\eqref{mhd}_4$ with respect to $t$ shows that
\begin{align}\label{lv4.12}
H_{tt}-H_t\cdot\nabla u-H\cdot\nabla u_t+u_t\cdot\nabla H+u\cdot\nabla H_t+H_t\divv u+H\divv u_t=\nu\Delta H_t.
\end{align}
Multiplying \eqref{lv4.12} by $H_t$ and integrating the resulting equality  over $B_R$ yield that
\begin{align}\label{lv4.13}
\frac12\frac{d}{dt}\int|H_t|^2dx+\nu\int|\nabla H_t|^2dx =&\int(H\cdot\nabla)u_t\cdot H_tdx+\int(u_t\cdot\nabla) H_t\cdot Hdx
 \notag \\
&+\int(H_t\cdot\nabla)u\cdot H_tdx-\frac{1}{2}\int |H_t|^2 \divv udx \triangleq \sum_{i=1}^4 S_i.
\end{align}
Integration by parts together with H{\"o}lder's inequality and \eqref{3.v2} leads to
\begin{align}\label{lvbo4.14}
S_1+S_2 & = -\int(H\cdot\nabla)H_t\cdot u_tdx+\int(u_t\cdot\nabla) H_t\cdot Hdx\notag \\
& \le 2\|\nabla H_t\|_{L^2} \||u_t||H|\|_{L^2} \notag \\
&\le \frac{\nu}{4}\|\nabla H_t\|_{L^2}^2
+C\|u_t\bar{x}^{-\frac{a}{4}}\|_{L^8}^2 \|H \bar{x}^{\frac{a}{2}}\|_{L^2}\|H\|_{L^4}
\notag \\
&\le \frac{\nu}{4}\|\nabla H_t\|_{L^2}^2+C\big(\|\sqrt{\rho}u_t\|_{L^2}^2+(1+\|\rho\|_{L^\infty}^2)\|\nabla u_t\|_{L^2}^2\big)\|H\bar{x}^{\frac{a}{2}}\|_{L^2}\|H\|_{L^4}\notag\\
& \le \frac{\nu}{4}\|\nabla H_t\|_{L^2}^2+C\big(\|H\bar{x}^{\frac{a}{2}}\|_{L^2}\|H\|_{L^4}(1+\|\rho\|_{L^\infty}^2)\big)\|\nabla u_t\|_{L^2}^2+\psi^\alpha.
\end{align}
By virtue of H{\"o}lder's inequality, \eqref{gjht}, and \eqref{2.4}, one has that
\begin{align}\label{lv4.14}
S_3+S_4 &  \leq \|H_t\|_{L^4}^2\|\nabla u\|_{L^2}
\leq  C \|H_t\|_{L^2}\|H_t\|_{H^1} \|\nabla u\|_{L^2}
\leq \frac{\nu}{4}\|\nabla H_t\|_{L^2}^2+C\psi^{\alpha}.
\end{align}
Inserting \eqref{lvbo4.14}--\eqref{lv4.14} into \eqref{lv4.13}, we get that
\begin{align}\label{ilv4.14}
& \frac{d}{dt}\|H_t\|_{L^2}^2+\nu\|\nabla H_t\|_{L^2}^2\le  C\big(\|H\bar{x}^{\frac{a}{2}}\|_{L^2}\|H\|_{L^4}(1+\|\rho\|_{L^\infty}^2)\big)\|\nabla u_t\|_{L^2}^2+\psi^\alpha.
\end{align}

2. From \eqref{gjht} and \eqref{3.22}, one has that
\begin{align*}
\|H_t(0)\|_{L^2}^2\le C E_0^{\alpha}\le C,
\end{align*}
which together with \eqref{ilv4.14}, Gronwall's inequality, \eqref{igj1}, \eqref{gj3}, \eqref{gj8}, and \eqref{gj6} yields that
\begin{align*}
&\sup_{0\le s\le t}\|H_s\|_{L^2}^2+\int_0^t\|\nabla H_s\|_{L^2}^2ds \notag\\
&\le C\sup_{0\le s\le t}\big(\|H\bar{x}^{\frac{a}{2}}\|_{L^2}\|H\|_{L^4}(1+\|\rho\|_{L^\infty}^2)\big)\int_0^t\|\nabla u_s\|_{L^2}^2ds+C\int_0^t\psi^\alpha ds\notag\\
&\le C\exp\left\{C\exp\left\{C\int_{0}^{t} \psi^{\alpha} ds\right\}\right\}.
\end{align*}
The proof of Lemma \ref{le--3} is finished.
 \hfill $\Box$

\begin{lemma}\label{lemma 3.8}
Let $(\rho,u,\theta,H)$ be as in Proposition \ref{pro} and $T_1$ be as in Lemma \ref{l3.01}. Then there exists a positive constant $\alpha>1$ such that, for all $t\in (0, T_1]$,
\begin{align}\label{wdgj}
\sup\limits_{0\le s\le t}\left(\|\sqrt{\rho}\theta\|_{L^2}^2+\|\nabla \theta\|_{L^2}^2+\|\sqrt{\rho}\theta_s\|_{L^2}^2\right) +\int_0^t\big(\|\sqrt{\rho}\theta_s\|_{L^2}^2+\|\nabla\theta_s\|_{L^2}^2\big)ds
\le C\exp\left\{C\int_{0}^{t} \psi^{\alpha} ds\right\}.
\end{align}
\end{lemma}
\textit{Proof.}
1. Choosing $b_1\le\frac{a}{2}$ in Lemma \ref{lemma2.5}, then for $0<b<\min\{b_1,1\}$, we have
$$\bar{x}^b\le C\big(1+|x|^{b_1}\big)<C\bar{x}^{\frac{a}{2}}.$$
Thus it follows from Lemma \ref{lemma2.5} that
\begin{align}\label{3.8.1}
& \int \Big[\frac{\mu}{2}|\nabla u+(\nabla u)^{tr}|^{2}+\lambda(\divv u)^2+\nu(\curl H)^{2}\Big]\bar{x}^b dx \notag \\
&\le C\|\nabla u\|_{L^2}^2+C\|\nabla H\|_{L^2}^2+C\int\left(\rho|\theta_t|+\rho|u||\nabla \theta|+\rho|\theta||\divv u|\right)|x|^{b_1}\notag \\
&\le C\psi^\alpha+C\|\sqrt{\rho} \bar{x}^{\frac{a}{2}}\|_{L^2\bigcap L^\infty}\left(\|\sqrt{\rho}\theta_t\|_{L^2}
+\|\sqrt{\rho} u\|_{L^2}\|\nabla\theta\|_{L^2}+\|\sqrt{\rho}\theta \|_{L^2}\|\nabla u\|_{L^2}\right)
\le C\psi^\alpha.
\end{align}
Multiplying \eqref{mhd}$_3$ by $\theta$ and integration by parts, one has that
\begin{align}\label{3.8.2}
&\frac{c_v}{2}\frac{d}{dt}\int\rho\theta^2dx+\kappa \int|\nabla \theta|^2dx\notag\\
&\le C\int\rho\theta^2|\divv u|dx+\int \Big[\frac{\mu}{2}|\nabla u+(\nabla u)^{tr}|^{2}+\lambda(\divv u)^2+\nu(\curl H)^{2}\Big]\theta dx.
\end{align}
By virtue of \eqref{zx} and \eqref{3.8.1}, we have
\begin{align}\label{bc1}
\int\rho\theta^2|\divv u|dx\le \|\rho\bar{x}^a\|_{L^q}\|\theta \bar{x}^{-\frac{a}{2}}\|_{\frac{4q}{q-2}}^2\|\nabla u\|_{L^2}\le C\psi^\alpha.
\end{align}
For simplicity, setting $Z\triangleq\frac{\mu}{2}|\nabla u+(\nabla u)^{tr}|^{2}+\lambda(\divv u)^2+\nu(\curl H)^{2}$, then we infer from \eqref{zx} and \eqref{3.8.1} that
\begin{align}\label{3.8.3}
\int Z\theta dx\le C\|\theta \bar{x}^{-\frac{b}{2}}\|_{L^6}
\|\sqrt{Z}\bar{x}^{\frac{b}{2}}\|_{L^2}
\big(\|\nabla u\|_{L^3}+\|\nabla H\|_{L^3}\big)
\le C\psi^\alpha,
\end{align}
due to $\sqrt{Z}\leq C(|\nabla u|+|\nabla H|)$.
Thus, by \eqref{bc1} and \eqref{3.8.3}, we obtain after integrating \eqref{3.8.2} over $[0,t]$ that
\begin{align}\label{3.8.5}
\sup\limits_{0\le s\le t} \|\sqrt{\rho} \theta\|_{L^2}^2
+\int_0^t\|\nabla\theta\|_{L^2}^2ds
\le C+C\int_0^t\psi^\alpha ds.
\end{align}

2. Multiplying \eqref{mhd}$_3$ by $\theta_t$ gives that
\begin{align}\label{3.8.6}
\frac{\kappa}{2}\frac{d}{dt}\|\nabla\theta\|_{L^2}^2
+c_v\|\sqrt{\rho}\theta_t\|_{L^2}^2
=-c_v\int (\rho u\cdot\nabla\theta) \theta_tdx
-R\int \rho\theta\theta_t\divv u+\int Z\theta_tdx.
\end{align}
By virtue of H{\"o}lder's inequality, \eqref{3.a2}, and \eqref{3.v2}, one has that
\begin{align}\label{3.76}
-c_v\int \rho u \cdot\nabla\theta\theta_t dx
& \le C\|\rho \bar{x}^a\|_{L^q}\|u \bar{x}^{-\frac{a}{2}}\|_{L^{\frac{4q}{q-2}}}\|\theta_t \bar{x}^{-\frac{a}{2}}\|_{L^{\frac{4q}{q-2}}}\|\nabla \theta\|_{L^2} \notag \\
& \le C\psi^\alpha\big(\|\sqrt{\rho}\theta_t\|_{L^2}
+(1+\|\rho\|_{L^\infty})\|\nabla\theta_t\|_{L^2}\big) \notag \\
& \le \frac{\kappa}{8}\|\nabla\theta_t\|_{L^2}^2+C\psi^\alpha,
\end{align}
and
\begin{align}
-R\int \rho\theta\theta_t\divv u\le C\|\sqrt{\rho}\theta\|_{L^\infty}\|\sqrt{\rho}\theta_t\|_{L^2}\|\nabla u\|_{L^2}\le C\psi^\alpha.
\end{align}
We deduce from H{\"o}lder's inequality, \eqref{3.v2}, and \eqref{3.8.1} that
\begin{align}\label{3.77}
\int Z\theta_t dx
\le C\|\theta_t \bar{x}^{-\frac{b}{2}}\|_{L^6}
\|\sqrt{Z}\bar{x}^\frac{b}{2}\|_{L^2}\big(\|\nabla u\|_{L^3}+\|\nabla H\|_{L^3}\big)
\le \frac{\kappa}{8}\|\nabla \theta_t\|_{L^2}^2+C\psi^\alpha.
\end{align}
Substituting \eqref{3.76}--\eqref{3.77} into \eqref{3.8.6} leads to
\begin{align}\label{3.8.7}
\frac{\kappa}{2}\frac{d}{dt}\|\nabla\theta\|_{L^2}^2
+c_v\|\sqrt{\rho}\theta_t\|_{L^2}^2
\leq\frac{\kappa}{4}\|\nabla \theta_t\|_{L^2}^2+C\psi^\alpha.
\end{align}

3. Differentiating \eqref{mhd}$_3$ with respect to $t$ and multiplying the resulting equation by $\theta_t$ yield that
\begin{align}\label{3.8.8}
\frac{c_v}{2}\frac{d}{dt}\|\sqrt{\rho}\theta_t\|_{L^2}^2+\kappa\|\nabla \theta_t\|_{L^2}^2&=-c_v\int\rho_t|\theta_t|^2 dx-c_v\int(\rho u)_t\cdot\nabla \theta\theta_t dx\notag\\
&\quad-\int P_t\divv u \theta_t-\int P\divv u_t\theta_t
+\int Z_t\theta_tdx \triangleq \sum_{i=1}^{5}L_i.
\end{align}
It follows from \eqref{mhd}$_1$, \eqref{3.53}, and integration by parts that
\begin{align*}
L_1
=-2c_v\int \rho u\cdot \nabla\theta_t \theta_tdx
\le\frac{\kappa}{16}\|\nabla \theta_t\|_{L^2}^2+C\|\sqrt{\rho} u\|_{L^\infty}^2\|\sqrt{\rho}\theta_t\|_{L^2}^2
\le\frac{\kappa}{16}\|\nabla \theta_t\|_{L^2}^2+C\psi^\alpha.
\end{align*}
In view of \eqref{mhd}$_1$, \eqref{3.v2}, and \eqref{3.a2}, we obtain from H{\"o}lder's inequality, \eqref{2.4}, and \eqref{2.5} that
\begin{align*}
L_2&=c_v\int (u \cdot\nabla \rho+\rho\divv u) \theta_t u\cdot\nabla \theta dx
-c_v\int\rho u_t\theta_t\nabla \theta dx\notag\\
&\le C\|\rho\bar{x}^a\|_{L^q}\|u \bar{x}^{-\frac{a}{3}}\|_{L^{\frac{6q}{q-2}}}^2\|\theta_t \bar{x}^{-\frac{a}{3}}\|_{L^{\frac{6q}{q-2}}}\|\nabla \theta\|_{L^2}+C\|\sqrt{\rho}u\|_{L^\infty}\|\nabla u\|_{L^4}\|\nabla \theta\|_{L^4}\|\sqrt{\rho}\theta_t\|_{L^2}\notag\\
&\quad+C\|\rho\|_{L^\infty}^{\frac14}\|\sqrt{\rho}\bar{x}^a\|_{L^\infty}^{\frac12}
\|\nabla \theta\|_{L^4}\|\theta_t\bar{x}^{-\frac{a}{2}}\|_{L^4}\|\sqrt{\rho}u_t\|_{L^2}
\notag\\
&\le\frac{\kappa}{16}\|\nabla\theta_t\|_{L^2}^2+C\psi^\alpha,
\end{align*}
and
\begin{align}
L_3 & =-R\int(\rho \theta)_t\divv u\theta_t dx\notag\\
&=R\int (u \cdot\nabla \rho+\rho \divv u)\theta\theta_t\divv udx-R\int\rho|\theta_t|^2\divv udx\notag\\
&\le\frac{\kappa}{16}\|\nabla\theta_t\|_{L^2}^2+C\psi^\alpha.
\end{align}
The Cauchy inequality together with \eqref{3.3.2} shows that
\begin{align}
L_4  =-R\int\rho \theta\divv u_t\theta_t dx
\le\|\nabla u_t\|_{L^2}^2+C\|\sqrt{\rho}\theta\|_{L^\infty}^2
\|\sqrt{\rho}\theta_t\|_{L^2}^2
\le\|\nabla u_t\|_{L^2}^2+C\psi^\alpha.
\end{align}
Direct calculation gives that
\begin{align*}
Z_t \leq C\sqrt{Z}\big(|\nabla u_t|+|\nabla H_t|\big),
\end{align*}
which combined with H{\"o}lder's inequality and \eqref{3.v2} ensures that
\begin{align}\label{sha7}
L_5 & \leq C\int|\theta_t|\sqrt{Z}\big(|\nabla u_t|+|\nabla H_t|\big)dx \notag \\
& \le C\|\theta_t\bar{x}^{-\frac{b}{4}}\|_{L^8}
\|Z^{\frac14}\bar{x}^\frac{b}{4}\|_{L^4}\|Z^{\frac14}\|_{L^8}
\||\nabla u_t|+|\nabla H_t|\|_{L^2}\notag\\
& \le C\|\theta_t\bar{x}^{-\frac{b}{4}}\|_{L^8}
\|\sqrt{Z}\bar{x}^\frac{b}{2}\|_{L^2}^{\frac12}
\big(\|\nabla u\|_{L^4}+\|\nabla H\|_{L^4}\big)^{\frac12}
\big(\|\nabla u_t\|_{L^2}+\|\nabla H_t\|_{L^2}\big) \notag\\
& \le C\exp\Big\{C\exp\Big\{C\int_{0}^{t} \psi^{\alpha} ds\Big\}\Big\}\Big(1+\|\sqrt{\rho}\theta_t\|_{L^2}^{\frac12}
+\|\nabla\theta\|_{L^2}^{\frac12}\Big)
\big(\|\nabla u_t\|_{L^2}^2+\|\nabla H_t\|_{L^2}^2\big) \notag \\
& \quad +\frac{\kappa}{16}\|\nabla \theta_t\|_{L^2}^2+C\psi^\alpha,
\end{align}
where in the last inequality we have used
\begin{align*}
\int Z\bar{x}^b dx
&\le C\|\sqrt{\rho} \bar{x}^{\frac{a}{2}}\|_{L^2\bigcap L^\infty}\big(\|\sqrt{\rho}\theta_t\|_{L^2}
+\|\sqrt{\rho} u\|_{L^2}\|\nabla\theta\|_{L^2}+\|\sqrt{\rho} \theta\|_{L^2}\|\nabla u\|_{L^2}\big)\notag \\
& \quad +C\|\nabla u\|_{L^2}^2+ C\|\nabla H\|_{L^2}^2 \notag \\
&\le  C\exp\left\{C\int_{0}^{t} \psi^{\alpha} ds\right\}\big(1+\|\sqrt{\rho}\theta_t\|_{L^2}+\|\nabla \theta\|_{L^2}\big),
\end{align*}
due to \eqref{3.8.1}, \eqref{gj3}, \eqref{wdgj}, and \eqref{gj8}.
Therefore, inserting the above estimates on $L_1$--$L_5$ into \eqref{3.8.8} and combining \eqref{3.8.7}, we find that
\begin{align}\label{z3.8.11}
& \frac{d}{dt}\big(c_v\|\sqrt{\rho}\theta_t\|_{L^2}^2
+\kappa\|\nabla\theta\|_{L^2}^2\big)
+\kappa\|\nabla \theta_t\|_{L^2}^2+c_v\|\sqrt{\rho}\theta_t\|_{L^2}^2
\notag \\
&\le C\exp\left\{C\int_{0}^{t} \psi^{\alpha} ds\right\}
\Big(1+\|\sqrt{\rho}\theta_t\|_{L^2}^{\frac12}
+\|\nabla \theta\|_{L^2}^{\frac12}\Big)
\big(\|\nabla u_t\|_{L^2}^2+\|\nabla H_t\|_{L^2}^2\big)+C\psi^\alpha.
\end{align}

4. It follows from \eqref{mhd}$_3$, \eqref{c.c}, and \eqref{3.53} that
\begin{align*}
\int \rho\theta_t^2(x,0)dx
& \le\lim_{t \rightarrow 0}\sup\Big(\int\rho^{-1}\Big[\kappa c_v^{-1}\Delta\theta
+\frac{\mu}{2c_v}|\nabla u+(\nabla u)^{tr}|^{2}\\
&\ \ \ +\lambda c_v^{-1}(\divv u)^{2}+\nu c_v^{-1}(\curl H)^{2}
-Rc_v^{-1}\rho\theta\divv u-\rho u\cdot\nabla\theta\Big]^2dx\Big) \\
& \le C\|g_2\|_{L^2}^2+C\|\rho u(0)\|_{L^\infty}^2\|\nabla\theta(0)\|_{L^2}^2 \\
& \leq C\|g_2\|_{L^2}^2+C E_0^\alpha \le C,
\end{align*}
which combined with \eqref{z3.8.11}  gives that,
for $t\in(0,T_1]$,
\begin{align}\label{3.8.11}
&\sup\limits_{0\le s\le t} \big(\|\sqrt{\rho} \theta_s\|_{L^2}^2+\|\nabla \theta\|_{L^2}^2\big)+\int_0^t\big(\|\sqrt{\rho}\theta_s\|_{L^2}^2+\|\nabla \theta_s\|_{L^2}^2\big)ds \notag \\
&\le C\exp\Big\{C\int_{0}^{t} \psi^{\alpha} ds\Big\}
\int_0^t \Big(1+\|\sqrt{\rho}\theta_s\|_{L^2}^{\frac12}
+\|\nabla \theta\|_{L^2}^{\frac12}\Big)
\big(\|\nabla u_s\|_{L^2}^2+\|\nabla H_s\|_{L^2}^2\big)ds \notag \\
& \quad +C\int_0^t\psi^\alpha ds+C.
\end{align}
By Young's inequality and \eqref{gj6}, it holds that
\begin{align}\label{3.80}
&C\exp\left\{C\int_{0}^{t} \psi^{\alpha} ds\right\}
\int_0^t\Big(1+\|\sqrt{\rho }\theta_s\|_{L^2}^{\frac12}+\|\nabla \theta\|_{L^2}^{\frac12}\Big)
\big(\|\nabla u_s\|_{L^2}^2+\|\nabla H_s\|_{L^2}^2\big)ds\notag \\
& \le\frac12\sup\limits_{0\le s\le t} \left(\|\sqrt{\rho} \theta_s\|_{L^2}^2+\|\nabla \theta\|_{L^2}^2\right) +C\exp\left\{C\int_{0}^{t} \psi^{\alpha} ds\right\}\Big(\int_0^t\big(\|\nabla u_s\|_{L^2}^2+\|\nabla H_s\|_{L^2}^2 \big)ds\Big)^{\frac43}\notag \\
& \le\frac{1}{2}\sup\limits_{0\le s\le t} \left(\|\sqrt{\rho} \theta_s\|_{L^2}^2+\|\nabla \theta\|_{L^2}^2\right)+C\exp\left\{C\int_{0}^{t} \psi^{\alpha} ds\right\}.
\end{align}
Thus, putting \eqref{3.80} into \eqref{3.8.11} leads to
\begin{align}\label{3.8.11''}
&\sup\limits_{0\le s\le t} \left(\|\sqrt{\rho}\theta_s\|_{L^2}^2+\|\nabla \theta\|_{L^2}^2\right)+\int_0^t\left(\|\sqrt{\rho}\theta_s\|_{L^2}^2+\|\nabla \theta_s\|_{L^2}^2\right)ds\le C\exp\left\{C\int_{0}^{t} \psi^{\alpha} ds\right\},
\end{align}
which along with \eqref{3.8.5} yields \eqref{wdgj}.
\hfill $\Box$

\begin{lemma}\label{xzlm}
Let $(\rho,u,\theta,H)$ be as in Proposition \ref{pro} and $T_1$ be as in Lemma \ref{l3.01}. Then there exists a positive constant $\alpha>1$ such that, for all $t\in (0, T_1]$,
\begin{align}\label{wdgj'}
&\sup\limits_{0\le s\le t}\left(\|\nabla^2 u\|_{L^2}^2+\|\nabla^2 H\|_{L^2}^2+\|\nabla^2 \theta\|_{L^2}^2\right) +\int_0^t\Big(\|\nabla^2\theta\|_{L^q}^2+\|\nabla^2\theta\|_{L^q}^\frac{q+1}{q}\Big)ds\notag\\
 &\le C\exp\left\{ C\exp\left\{C\int_{0}^{t} \psi^{\alpha} ds\right\}\right\}.
\end{align}
\end{lemma}
\textit{Proof.}
1. It deduces from $\eqref{mhd}_4$, the standard $L^2$-estimate of elliptic equations, \eqref{3.a2}, H\"older's inequality, \eqref{3.1.0}, and Gagliardo-Nirenberg inequality that
\begin{align}\label{AMSS11}
\|\nabla^2H\|^2_{L^2}
&\leq C\|H_t\|^2_{L^2}+C\||u||\nabla H|\|_{L^2}^2+C\||H||\nabla u|\|^2_{L^2} \notag \\
& \leq C\|H_t\|^2_{L^2}+C\|u \bar{x}^{-\frac{a}{4}}\|_{L^8}^2\|\nabla H \bar{x}^{\frac{a}{2}}\|_{L^2}\|\nabla H\|_{L^4}
+C\|H\|_{L^4}^2\|\nabla u\|_{L^4}^2 \notag \\
& \leq C\|H_t\|^2_{L^2}+C\|\nabla H \bar{x}^{\frac{a}{2}}\|_{L^2}^2
+C\big(1+(1+\|\rho\|_{L^\infty}^4)\|\nabla u\|_{L^2}^4\big)\|\nabla H\|_{L^2}\|\nabla H\|_{H^1}\notag\\
&\quad+C\|H\|_{L^4}^2\|\nabla u\|_{L^2}\|\nabla u\|_{H^1} \notag \\
& \leq C\|H_t\|^2_{L^2}+C\|\nabla H \bar{x}^{\frac{a}{2}}\|_{L^2}^2
+\frac14\|\nabla^2 H\|_{L^2}^2+\frac14\|\nabla^2 u\|_{L^2}^2+C\|H\|_{L^4}^8 \notag \\
& \quad +C\big(1+(1+\|\rho \bar{x}^a\|_{W^{1,q}\cap H^1}^8)\|\nabla u\|^8_{L^2}\big)\big(1+\|\nabla H\|^2_{L^2}\big).
\end{align}
It follows from \eqref{lv3.60} with $p=2$, \eqref{3.v2}, \eqref{3.3.1}, \eqref{3.a2}, and Gagliardo-Nirenberg inequality that
\begin{align}\label{bc2}
\|\nabla^2u\|_{L^2}^2
\leq& \|\rho u_t\|_{L^2}^2+C\|\rho u\cdot\nabla u\|_{L^2}^2
+C\||H||\nabla H|\|_{L^2}^2+C\|\nabla P\|_{L^2}^2
\notag \\
\leq& C \|\rho\bar{x}^a\|_{{W^{1,q}\cap H^1}}\|\sqrt{\rho}u_t\|_{L^2}^2
+C\|\rho u\|_{L^4}^2\|\nabla u\|_{L^4}^2+C\|H\|_{L^4}^2\|\nabla H\|_{L^4}^2 +C\|\nabla \rho \theta\|_{L^2}^2+C\|\rho\nabla\theta\|_{L^2}^2\notag \\
\leq& C \|\rho\bar{x}^a\|_{{W^{1,q}\cap H^1}}\big(\|\sqrt{\rho}u_t\|_{L^2}^2
+\big(1+(1+\|\rho\|_{L^\infty}^2)\|\nabla u\|_{L^2}^2\big)\|\nabla u\|_{L^2}\|\nabla u\|_{H^1}\big)\notag \\
&+C\|H\|_{L^4}^2\|\nabla H\|_{L^2}\|\nabla H\|_{H^1}+C\|\rho\bar{x}^a\|_{W^{1,q}\cap H^1}^2\Big(\|\theta\bar{x}^{-a}\|_{L^{\frac{2q}{q-2}}}^2
+\|\nabla\theta\|_{L^2}^2\Big)\notag\\
\leq& C \|\sqrt{\rho}u_t\|_{L^2}^4
+C\big(1+(1+\|\rho\bar{x}^a\|_{W^{1,q}\cap H^1}^{8})\|\nabla u\|_{L^2}^{12}\big)
+\frac14\|\nabla^2u\|_{L^2}^2+\frac14\|\nabla^2H\|_{L^2}^2\notag\\
&+C\|\rho\bar{x}^a\|_{W^{1,q}\cap H^1}^2\big(1+\|\sqrt{\rho}\theta\|_{L^2}^2+\big(1+\|\rho\bar{x}^a\|_{W^{1,q}\cap H^1}^{2}\big)\|\nabla\theta\|_{L^2}^2\big)\notag\\
&+C\|\nabla H\|_{L^2}^4+C\|H\|_{L^4}^8+C\|\rho\bar{x}^a\|_{W^{1,q}\cap H^1}^4.
\end{align}
This along with \eqref{gj3}, \eqref{gj6}, \eqref{AMSS11}, \eqref{gj10'}, \eqref{gj8}, and \eqref{igj10'} yields that
\begin{align}\label{iAMSS12}
\sup_{0\le s\le t}\big(\|\nabla^2u\|_{L^2}^2+\|\nabla H\|_{L^2}^2
\big)
 \leq C\exp\left\{C\exp\left\{C\int_{0}^{t}\psi^{\alpha}ds\right\}\right\}.
\end{align}
We deduce from \eqref{mhd}$_3$, the standard $L^2$-estimate of elliptic equations, \eqref{3.3.1}, \eqref{3.a2}, and \eqref{2.4} that
\begin{align}\label{3.8.12}
 \|\nabla^2\theta\|_{L^2}^2
&\le C \big(\|\rho\theta_t\|_{L^2}^2+ \|\rho u\cdot\nabla\theta\|_{L^2}^2+C\|\rho\theta\divv u\|_{L^2}^2
+ \|\nabla u\|_{L^4}^4+ \|\nabla H\|_{L^4}^4\big) \notag \\
&\le C \big(\|\rho\bar{x}^a\|_{{W^{1,q}\cap H^1}}\|\sqrt{\rho}\theta_t\|_{L^2}^2
+\|\rho u\|_{L^4}^2\|\nabla\theta\|_{L^4}^2
+\|\rho \theta\|_{L^4}^2\|\nabla u\|_{L^4}^2\notag\\
&\quad+\|\nabla u\|_{L^2}^2\|\nabla u\|_{H^1}^2
+\|\nabla H\|_{L^2}^2\|\nabla H\|_{H^1}^2\big) \notag \\
&\le C\|\rho\bar{x}^a\|_{{W^{1,q}\cap H^1}}\big(\|\sqrt{\rho}\theta_t\|_{L^2}^2
+\big(1+(1+\|\rho\bar{x}^a\|_{W^{1,q}\cap H^1}^{2})\|\nabla u\|_{L^2}^2\big)\|\nabla\theta\|_{L^2}\|\nabla\theta\|_{H^1}\big)\notag\\
&\quad+C\|\rho\bar{x}^a\|_{{W^{1,q}\cap H^1}}\big(\|\sqrt{\rho}\theta\|_{L^2}^2+(1+\|\rho\bar{x}^a\|_{W^{1,q}\cap H^1}^{2})\|\nabla \theta\|_{L^2}^2\big)\|\nabla u\|_{L^2}\|\nabla u\|_{H^1}
\notag \\
&\quad +C\big(\|\nabla u\|_{L^{2}}^4+ \|\nabla^2 u\|_{L^2}^4+\|\nabla H\|_{L^2}^4+ \|\nabla^2 H\|_{L^2}^4\big) \notag \\
&\le C\|\sqrt{\rho}\theta_t\|_{L^2}^4
+\frac12\|\nabla^2\theta\|_{L^2}^2+C\big(1+(1+\|\rho\bar{x}^a\|_{W^{1,q}\cap H^1}^{8})\|\nabla u\|_{L^2}^8\big)
\|\nabla\theta\|_{L^2}^4\notag\\
&\quad+\big(1+\|\sqrt{\rho}\theta\|_{L^2}^4+(1+\|\rho\bar{x}^a\|_{W^{1,q}\cap H^1}^{2})\|\nabla\theta\|_{L^2}^4\big)\|\rho\bar{x}^a\|_{{W^{1,q}\cap H^1}}^4 \notag\\
&\quad +C\big(\|\nabla u\|_{L^{2}}^4+ \|\nabla^2 u\|_{L^2}^4+\|\nabla H\|_{L^2}^4+ \|\nabla^2 H\|_{L^2}^4+\|\rho\bar{x}^a\|_{{W^{1,q}\cap H^1}}^4\big),
\end{align}
which together with \eqref{gj3}, \eqref{gj8}, \eqref{wdgj}, and \eqref{iAMSS12} leads to
\begin{align}\label{iiAMSS12}
\sup_{0\le s\le t}\|\nabla^2\theta\|_{L^2}^2
 \leq C\exp\left\{C\exp\left\{C\int_{0}^{t}\psi^{\alpha}ds\right\}\right\}.
\end{align}
The standard $L^q$-estimate of elliptic equations together with \eqref{mhd}$_3$, \eqref{3.2}, H{\"o}lder's inequality, \eqref{local1}, and \eqref{2.4} yields that
\begin{align}\label{cz}
 \|\nabla^2\theta\|_{L^q}
& \le C \left(\|\rho \theta_t\|_{L^q}+\|\rho u\cdot\nabla\theta\|_{L^q}+\|\rho\theta\divv u\|_{L^q}+ \|\nabla u\|_{L^{2q}}^2+\|\nabla H\|_{L^{2q}}^2\right) \notag \\
& \le C\|\rho\|_{L^\infty}^{\frac{1}{2}}\|\sqrt{\rho}\theta_t\|_{L^2}^{\frac{2(q-1)}{q^2-2}}
\|\sqrt{\rho}\theta_t\|_{L^{q^2}}^{\frac{q^2-2q}{q^2-2}}
+C\|\rho u\|_{L^{2q}}\|\nabla\theta\|_{L^{2q}}+C\|\rho \theta\|_{L^{2q}}\|\nabla u\|_{L^{2q}} \notag \\
& \quad +C\|\nabla u\|_{L^2}^{\frac2q}\|\nabla u\|_{H^1}^{\frac{2q-2}{q}}
+C\|\nabla H\|_{L^2}^{\frac2q}\|\nabla H\|_{H^1}^{\frac{2q-2}{q}}
\notag \\
& \le C\psi^\alpha\|\sqrt{\rho}\theta_t\|_{L^2}^{\frac{2(q-1)}{q^2-2}}
\big(\|\sqrt{\rho}\theta_t\|_{L^2}+(1+\|\rho\|_{L^\infty})\|\nabla\theta_t\|_{L^2}\big)^{\frac{q^2-2q}{q^2-2}}\notag\\
&\quad+C\psi^\alpha\big(\|\sqrt{\rho}u\|_{L^2}+(1+\|\rho\|_{L^\infty})\|\nabla u\|_{L^2}\big)
\|\nabla \theta\|_{L^2}^{\frac1q}\|\nabla\theta\|_{H^1}^{\frac{q-1}{q}}\notag \\
& \quad +C\psi^\alpha\big(\|\sqrt{\rho}\theta\|_{L^2}+(1+\|\rho\|_{L^\infty})\|\nabla \theta\|_{L^2}\big)\|\nabla u\|_{L^2}^{\frac1q}\|\nabla u\|_{H^1}^{\frac{q-1}{q}}\notag\\
&\quad+C\|\nabla u\|_{L^2}^{\frac2q}\|\nabla u\|_{H^1}^{\frac{2q-2}{q}}
+C\|\nabla H\|_{L^2}^{\frac2q}\|\nabla H\|_{H^1}^{\frac{2q-2}{q}} \notag \\
& \le C\psi^\alpha
\|\nabla\theta_t\|_{L^2}^{\frac{q^2-2q}{q^2-2}}
+C\psi^\alpha.
\end{align}
Consequently, similarly to \eqref{olv3.61} and \eqref{olv3.62}, we infer from \eqref{cz}, Young's inequality, and \eqref{3.8.11''} that
\begin{align} \label{3.8.16}
\int_0^t\Big(\|\nabla^2\theta\|_{L^q}^{\frac{q+1}{q}}
+\|\nabla^2\theta\|_{L^q}^2\Big)ds
\leq C\int_0^t\big(\psi^\alpha+\|\nabla\theta_s\|_{L^2}^2\big)ds
\le C\exp\left\{C\int_{0}^{t} \psi^{\alpha} ds\right\},
\end{align}
which together with \eqref{iAMSS12} and \eqref{iiAMSS12} gives \eqref{wdgj'}. The proof of Lemma \ref{xzlm} is finished.   \hfill $\Box$

Now, Proposition \ref{pro} is a direct consequence of  Lemmas \ref{l3.01}--\ref{xzlm}.

\textit{Proof of Proposition \ref{pro}.}
It follows from \eqref{igj1}, \eqref{gj3}, \eqref{gj6}, \eqref{gj10'}, \eqref{gj8}, \eqref{wdgj}, and \eqref{wdgj'} that
\begin{align*}
\psi(t)\le C\exp\left\{C\exp\left\{C\int_{0}^{t} \psi^{\alpha} ds\right\}\right\}.
\end{align*}
Standard arguments yield  that for $M\triangleq Ce^{Ce}$ and $T_0\triangleq \min\{T_1,(CM^\alpha)^{-1}\}$,
\begin{align*}
\sup\limits_{0\le t\le T_0}\psi(t)\le M,
\end{align*}
which together with \eqref{igj1}, \eqref{gj3}, \eqref{gj6}, \eqref{gj10'}, \eqref{gj8}, \eqref{gj9'}, \eqref{wdgj}, and \eqref{wdgj'} gives \eqref{o1}. The proof of Proposition \ref{pro} is complete.   \hfill $\Box$

\section{Proof of Theorem \ref{t1}}\label{sec4}

With the \textit{a priori} estimates in Section \ref{sec3} at hand, we are now in a position to prove Theorem \ref{t1}.

\textbf{Step 1. Local existence of strong solutions.}
Let $(\rho_{0},u_{0},\theta_{0},H_{0})$ be as in Theorem \ref{t1}. Without loss of generality, we assume that the initial density $\rho_0$ satisfies
\begin{align*}
\int_{\mathbb{R}^2} \rho_0dx=1,
\end{align*}
which indicates that there exists a positive constant $N_0$ such that
\begin{align}\label{oi3.8}
\int_{B_{N_0}} \rho_0 dx\ge \frac34\int_{\mathbb{R}^2}\rho_0dx=\frac34.
\end{align}
We construct
$\rho_{0}^{R}=\hat\rho_{0}^{R}+R^{-1}e^{-|x|^2} $, where $0\le\hat\rho_{0}^{R}\in  C^\infty_0(\mathbb{R}^2)$ satisfies
\begin{align}\label{bci0}
\begin{cases}
\int_{B_{N_0}}\hat\rho^R_0dx\ge \frac12,\\
\bar x^a \hat\rho_{0}^{R}\rightarrow \bar x^a \rho_{0}\quad {\rm in}\,\, L^1(\mathbb{R}^2)\cap H^{1}(\mathbb{R}^2)\cap W^{1,q}(\mathbb{R}^2)\ \ {\rm as}\ \ R\rightarrow\infty.
\end{cases}
\end{align}
Noting that $H_0\bar{x}^{\frac{a}{2}}\in H^1(\mathbb{R}^2)$ and
$\nabla H_0\in H^1(\mathbb{R}^2)$, we choose $H_0^R\in \{w\in C^\infty_0(B_R)~|~\divv w=0\}$ satisfying
\begin{align}\label{lv6.1}
H_0^R\bar x^{\frac{a}{2}} \rightarrow  H_0\bar x^{\frac{a}{2}},\
 \nabla H_{0}^{R}\rightarrow \nabla H_{0}\ \
{\rm in}\ \ H^1(\mathbb{R}^2)\ \ {\rm as}\ \ R\rightarrow\infty.
\end{align}
Since $\nabla u_0\in H^1(\mathbb{R}^2),$  we select $v^R_i\in C^\infty_0(B_R)~(i=1,2)$ such that for $i=1,2,$
\begin{align}\label{bci3}
\lim\limits_{R\rightarrow \infty}\|v^R_i-\partial_iu_0\|_{H^1(\mathbb{R}^2)}=0.
\end{align}

Consider the unique smooth solution $u_0^R$ to the elliptic problem
\begin{align}\label{bbi2}
\begin{cases}
-\mu\Delta u_{0}^{R}-(\mu+\lambda) \nabla \divv u_{0}^{R}+\rho_0^R u_0^R+\nabla P^{R}_0 \\
= H_{0}^{R}\cdot \nabla H_{0}^{R}-\frac{1}{2}\nabla| H_{0}^{R}|^2+\sqrt{\rho_{0}^{R}} h_1^R-  \partial_iv^R_i,& {\rm in} \,\,  B_{R},\\
u_{0}^{R} =0,\,\, \,& {\rm on}~\partial B_{R},
\end{cases}
\end{align}
where $h_1^R=(\sqrt{\rho_0}u_0)*j_{\frac1R}$ with $j_\delta$ being the standard mollifying kernel of width $\delta$.
Extending $u_{0}^{R} $ to $\mathbb{R}^2$ by defining $0$ outside $B_{R}$ and denoting it by $\tilde{u}_{0}^{R}$, we claim that, up to the extraction of subsequences,
\begin{align}\label{3.7i4}
\lim\limits_{R\rightarrow \infty}
\left(\big\|\nabla\tilde{u}_{0}^{R}-\nabla u_0\big\|_{H^1(\mathbb{R}^2)}
+\Big\|\sqrt{\rho_0^R}\tilde{u}_{0}^{R}-\sqrt{\rho_0}u_0\Big\|_{L^2(\mathbb{R}^2)}\right)=0.
\end{align}
Indeed, it is not hard to find that $\tilde{u}_0^R$ is also a solution of \eqref{bbi2} in $\mathbb{R}^2$. Multiplying \eqref{bbi2} by $\tilde{u}_0^R$ and integrating the resulting equation over $\mathbb{R}^2$ lead to
\begin{align*}
&\big\|\sqrt{\rho_0^R}\tilde{u}_0^R\big\|_{L^2(\mathbb{R}^2)}^2
+\mu\|\nabla \tilde{u}_0^R\|_{L^2(\mathbb{R}^2)}^2+(\mu+\lambda)\|\divv \tilde{u}_0^R\|_{L^2(\mathbb{R}^2)}^2\\
&\le \big\|\sqrt{\rho_0^R}\tilde{u}_0^R\big\|_{L^2(B_R)}\|h_1^R\|_{L^2(B_R)}+C\|P_0^R\|_{L^2(B_R)}\|\nabla\tilde{u}_0^R\|_{L^2(B_R)}\\
&\quad +\|v_i^R\|_{L^2(B_R)}\|\partial_i\tilde{u}^R_0\|_{L^2(B_R)}+C\|H_0^R\|_{L^4(B_R)}^2\|\nabla\tilde{u}_0^R\|_{L^2(B_R)}+C\||H_0^R||\nabla H_0^R||\tilde{u}_0^R|\|_{L^1(B_R)}\\
&\le \frac12\big\|\sqrt{\rho_0^R}\tilde{u}_0^R\big\|_{L^2(B_R)}^2+\frac{\mu}2\|\nabla \tilde{u}_0^R\|_{L^2(B_R)}^2+\frac12\|h_1^R\|_{L^2(B_R)}^2
+C\|\rho_0^R\|_{L^\infty(B_R)}\big\|\sqrt{\rho_0^R}\theta_0^R \big\|_{L^2(B_R)}^2
\\
&\quad+C\|v_i^R\|_{L^2(B_R)}^2
+C\|H_0^R\|_{L^4(B_R)}^4+C\|H_0^R\|_{L^4(B_R)}\|\bar{x}^{-\frac{a}{2}}
\tilde{u}_0^R\|_{L^4(B_R)}\|\bar{x}^{\frac{a}{2}}\nabla H_0^R\|_{L^2(B_R)}\\
&\le \frac12\big\|\sqrt{\rho_0^R}\tilde{u}_0^R\big\|_{L^2(B_R)}^2+\frac{\mu}2\|\nabla \tilde{u}_0^R\|_{L^2(B_R)}^2+C+C\big(\big\|\sqrt{\hat{\rho}_0^R}\tilde{u}_0^R\big\|_{L^2(B_R)}+(1+\|\hat{\rho}_0^R\|_{L^\infty})\|\nabla \tilde{u}_0^R\|_{L^2(B_R)}\\
&\le \frac12\big\|\sqrt{\rho_0^R}\tilde{u}_0^R\big\|_{L^2(B_R)}^2+\frac{\mu}2\|\nabla \tilde{u}_0^R\|_{L^2(B_R)}^2+C,
\end{align*}
owing to \eqref{bci0}, \eqref{lv6.1}, \eqref{o1}, and \eqref{3.v2}, which implies that
\begin{align}\label{2.i9-4}
\big\|\sqrt{\rho_0^R}\tilde{u}_0^R\big\|_{L^2(\mathbb{R}^2)}^2
+\|\nabla \tilde{u}_0^R\|_{L^2(\mathbb{R}^2)}^2\le C.
\end{align}
This together with \eqref{bci0} yields that there exist a subsequence $R_j\rightarrow \infty$ and a function $\tilde{u}_0\in\{\tilde{u}_0\in H^1_{\rm loc}(\mathbb{R}^2)|\sqrt{\rho_0}\tilde{u}_0\in L^2(\mathbb{R}^2), \nabla \tilde{u}_0\in L^2(\mathbb{R}^2)\}$ such that
\begin{align}\label{bci9}
\sqrt{\rho^{R_j}_0}\tilde{u}^{R_j}_0 \rightharpoonup\sqrt{\rho_0}\tilde{u}_0,\
\nabla\tilde{u}_0^{R_j}\rightharpoonup \nabla\tilde{u}_0\ \mbox{ weakly in } L^2(\mathbb{R}^2).
\end{align}
We claim that
\begin{align}\label{bai1}
\tilde{u}_0=u_0.
\end{align}
In fact, subtracting \eqref{c.c} from \eqref{bbi2} gives
\begin{align}\label{bai3}
-\mu&\Delta\left(\tilde{u}_0^{R_j}-u_0\right)-(\mu+\lambda)\nabla\divv\left(\tilde{u}_0^{R_j}-u_0\right)+\nabla\left(\tilde{P}_0^{R_j}-P_0\right)\notag\\
&=H_0^{R_j}\cdot\nabla H_0^{R_j}-H_0\cdot\nabla H_0-\frac{1}{2}\nabla\left(|H_0^{R_j}|^2-|H_0|^2\right)\notag\\
&\quad +\left(\sqrt{\rho_0^{R_j}}g_1*j_{1/{R_j}}-\sqrt{\rho_0}g_1\right)-\sqrt{\rho_0^{R_j}}
\left(\sqrt{\rho_0^{R_j}}\tilde{u}_0^{R_j}-\sqrt{\rho_0}u_0*j_{1/{R_j}}\right).
\end{align}
Multiplying \eqref{bbi2} by a test function $\pi\in C_0^\infty(\mathbb{R}^2)$ , it holds that
\begin{align}\label{4.11}
&\mu\int_{\mathbb{R}^2}\partial_i\left(\tilde{u}_0^{R_j}-v_i^{R_j}\right)
\cdot\partial_i\pi dx+(\mu+\lambda)\int_{\mathbb{R}^2} \divv(\tilde{u}^{R_j}_0-u_0)\divv\pi dx+\int_{\mathbb{R}^2} \sqrt{\rho^{R_j}_0}\Big(\sqrt{\rho^{R_j}_0}\tilde{u}^{R_j}_0-h_1^{R_j}\Big)\cdot\pi dx\notag\\
&=\int_{\mathbb{R}^2}\left(P_0^{R_j}-P_0\right)\divv \pi dx+\int\left(\sqrt{\rho_0^{R_j}}g_1*j_{\frac{1}{R_j}}-\sqrt{\rho_0}g_1\right)\pi dx+\int\left(H_0^{R_j}-H_0\right)\cdot\nabla H_0^{R_j} \pi dx\notag\\
&\quad+\int H_0^{R_j}\cdot\nabla\left(H_0^{R_j}-H_0\right) \pi dx+\frac{1}{2}\int\Big(|H_0^{R_j}|^2-|H_0|^2\Big)\divv \pi dx.
\end{align}
Then, letting $R_j\rightarrow \infty$ in \eqref{4.11}, it follows from \eqref{bci0}, \eqref{lv6.1}, \eqref{bci3}, and \eqref{bci9} that
\begin{align}\label{ibai2}
\int_{\mathbb{R}^2}\partial_i(\tilde{u}_0-u_0)\cdot\partial_i\pi dx+\int_{\mathbb{R}^2} \rho_0(\tilde{u}_0-u_0)\cdot\pi dx=0,
\end{align}
which implies \eqref{bai1} for the arbitrary $\pi$. Furthermore, multiplying \eqref{bbi2} by $\tilde{u}_0^{R_j}$ and integrating the resulting equation over $\mathbb{R}^2$, by the same arguments as \eqref{ibai2}, we have
\begin{align*}
&\lim\limits_{R_j\rightarrow\infty}
\left(\mu\Big\|\nabla\tilde{u}_0^{R_j}\Big\|_{L^2(\mathbb{R}^2)}^2
+(\mu+\lambda)\Big\|\divv\tilde{u}_0^{R_j}\Big\|_{L^2(\mathbb{R}^2)}^2+\Big\|\sqrt{\rho_0^{R_j}}\tilde{u}_0^{R_j}\Big\|_{L^2(\mathbb{R}^2)}^2\right)\\
&=\Big(\mu\|\nabla u_0\|_{L^2(\mathbb{R}^2)}^2+(\mu+\lambda)\|\divv u_0\|_{L^2(\mathbb{R}^2)}^2+\|\sqrt{\rho_0}u_0\|_{L^2(\mathbb{R}^2)}^2\Big),
\end{align*}
which combined with \eqref{bci9} leads to
\begin{align*}
\lim\limits_{R_j\rightarrow\infty}\big\|\nabla\tilde{u}_0^{R_j}\big\|_{L^2(\mathbb{R}^2)}^2
=\|\nabla\tilde{u}_0\|_{L^2(\mathbb{R}^2)}^2,\
\lim\limits_{R_j\rightarrow\infty}
\Big\|\sqrt{\rho_0^{R_j}}\tilde{u}_0^{R_j}\Big\|_{L^2(\mathbb{R}^2)}^2
=\|\sqrt{\rho_0}\tilde{u}_0\|_{L^2(\mathbb{R}^2)}^2.
\end{align*}
This along with \eqref{bai1} and \eqref{bci9} guarantees that
\begin{align}\label{4.13}
\lim\limits_{R_j\rightarrow \infty}
\left(\big\|\nabla\tilde{u}_{0}^{R_j}-\nabla u_0\big\|_{L^2(\mathbb{R}^2)}
+\Big\|\sqrt{\rho_0^{R_j}}\tilde{u}_{0}^{R_j}-\sqrt{\rho_0}u_0\Big\|_{L^2(\mathbb{R}^2)}\right)=0.
\end{align}
Moreover, if we differentiate \eqref{bbi2} and then multiply the resultant equality by $\nabla\tilde{u}_0^{R}$, it is not hard to infer from \eqref{c.c}, \eqref{bci9}, and \eqref{bai1} that, up to the extraction of subsequences,
\begin{align}\label{4.14}
\lim\limits_{R_j\rightarrow \infty}
\big\|\nabla^2\tilde{u}_{0}^{R_j}-\nabla^2u_0\big\|_{L^2(\mathbb{R}^2)}=0,
\end{align}
which together with \eqref{4.13} implies \eqref{3.7i4}.

Next, we consider the following elliptic equation on $\theta_0^R$:
\begin{align}\label{bbi3}
\begin{cases}
-\kappa\Delta \theta_{0}^{R}+\rho_0^R \theta_0^R
=\frac{\mu}{2}|\nabla u_0^R+(\nabla u_0^R)^{tr}|^2+\lambda(\divv u_0)^2
+\nu(\curl H_0^R)^2+\sqrt{\rho_0^R}h_2^R ,& {\rm in} \,\,  B_{R},\\
\theta_{0}^{R} =0,\,\, \,& {\rm on}~\partial B_{R},
\end{cases}
\end{align}
where $h_2^R=(\sqrt{\rho_0}\theta_0+g_2)*j_{\frac1R}$. Extending $\theta_{0}^{R} $ to $\mathbb{R}^2$ by defining $0$ outside $B_{R}$ and denoting it by $\tilde{\theta}_{0}^{R}$.
Similarly to \eqref{4.13} and \eqref{4.14}, we can also deduce that (see \cite{L2021})
\begin{align}\label{4.7i4}
\lim\limits_{R\rightarrow \infty}
\left(\big\|\nabla\tilde{\theta}_{0}^{R}-\nabla\theta_0\big\|_{H^1(\mathbb{R}^2)}
+\Big\|\sqrt{\rho_0^R}\tilde{\theta}_{0}^{R}-\sqrt{\rho_0}\theta_0 \Big\|_{L^2(\mathbb{R}^2)}\right)=0.
\end{align}
Hence, by virtue of Lemma \ref{th0}, the initial-boundary-value problem  \eqref{mhd} and \eqref{2.1} with the initial data $(\rho_0^R,u_0^R,\theta_0^R,H_0^R)$ has a unique strong solution $(\rho^{R},u^{R},\theta^R,H^{R})$ on $B_{R}\times (0,T_R]$. Moreover, Proposition \ref{pro} shows that there exists a $T_0$ independent of $R$ such that \eqref{o1} holds for $(\rho^{R},u^{R},\theta^R,H^{R})$.

For simplicity, in what follows, we write
\begin{align*}
L^p=L^p(\mathbb{R}^2),\ W^{k,p}=W^{k,p}(\mathbb{R}^2),\ H^k=W^{k,2}.
\end{align*}
Extending $(\rho^{R},u^{R},\theta^R,H^{R})$ by zero on $\mathbb{R}^2\setminus B_{R}$ and denoting it by
\begin{align*}
\big(\tilde{\rho}^R\triangleq \varphi_R\rho^R,  \tilde{u}^R,\tilde{\theta}^R,  \tilde{H}^R\big)
\end{align*}
with $\varphi_R$ satisfying \eqref{vp1}. From \eqref{o1}, we have
\begin{align}\label{kq1}
& \sup\limits_{0\le t\le T_0}
\big(\|\sqrt{\tilde{\rho}^R }\tilde{u}^R\|_{L^2}+\|\sqrt{\tilde{\rho}^R }\tilde{\theta}^R\|_{L^2}+\|\nabla\tilde{u}^R\|_{H^1}+\|\nabla \tilde{\theta}^R\|_{H^1}+\|\tilde{H}^R\|_{H^2}
+\|\tilde{H}^R\bar x^{\frac{a}{2}}\|_{H^1}\big) \notag \\
& \le \sup\limits_{0\le t\le T_0}\big(\|\sqrt{\rho^R}u^R\|_{L^2(B_R)}+\|\sqrt{\rho^R}\theta^R\|_{L^2(B_R)}
+\|\nabla u^R\|_{H^1(B_R)}+\|\nabla\theta^R\|_{H^1(B_R)}\notag\\
&\quad+\|H^R\|_{H^2(B_R)}+\|H^R\bar x^{\frac{a}{2}}\|_{H^1(B_R)}\big)
\le C,
\end{align}
and
\begin{align*}
\sup\limits_{0\le t\le T_0}\|\tilde{\rho}^R\bar x^a\|_{L^1}\le C.
\end{align*}
Similarly, it follows from \eqref{o1} that, for $q$ as in Theorem \ref{t1},
\begin{align}
& \sup\limits_{0\le t\le T_0}\Big(\big\|\sqrt{\tilde{\rho}^R}\tilde{u}^R_t\big\|_{L^2}
+\big\|\sqrt{\tilde{\rho}^R}  \tilde{\theta}^R_t\big\|_{L^2}
+\big\|\tilde{H}^R_t\big\|_{L^2}\Big) \notag \\
& \quad +\int_0^{T_0}\Big(\big\|\sqrt{\tilde{\rho}^R}\tilde{u}^R_t\big\|_{L^2}^2
+\big\|\sqrt{\tilde{\rho}^R}\tilde{\theta}^R_t\big\|_{L^2}^2
+\big\|\tilde{H}^R_t\big\|_{L^2}^2
\Big)dt \notag \\
& \quad +\int_0^{T_0}\Big(\big\|\nabla^2\tilde{u}^R\big\|_{L^q}^{\frac{q+1}{q}}
+\big\|\nabla^2\tilde{\theta}^R\big\|_{L^q}^{\frac{q+1}{q}}
+\big\|\nabla^2\tilde{u}^R\big\|_{L^q}^{2}
+\big\|\nabla^2\tilde{\theta}^R\big\|_{L^q }^{2}
+\big\|\nabla\tilde{u}^R_t\big\|_{L^2}^2\Big) dt \notag \\
& \quad +\int_0^{T_0}\Big(\big\|\nabla^2\tilde{H}^R\bar x^{\frac{a}{2}}\big\|_{L^2}^2
+\big\|\nabla\tilde{H}^R\bar x^{\frac{a}{2}}\big\|_{L^2}^2
+\big\|\nabla\tilde{\theta}^R_t\big\|_{L^2}^2+\big\|\nabla\tilde{H}^R_t\big\|_{L^2}^2\Big)dt
\le C.
\end{align}
Next, for $p\in[2,q]$, we obtain from \eqref{o1} and \eqref{gj8} that
\begin{align}
\sup\limits_{0\le t\le T_0}\big\|\nabla(\tilde{\rho}^R\bar x^a)\big\|_{L^p }
&\le C\sup\limits_{0\le t\le T_0}\Big(\big\|\nabla(\rho^R\bar x^a)\big\|_{L^p(B_R}
+R^{-1}\big\|\rho^R\bar x^a\big\|_{L^p(B_R)}\Big) \notag \\
& \le C\sup\limits_{0\le t\le T_0}\big\|\rho^R\bar x^a\big\|_{H^1(B_R)\cap W^{1,p}(B_R)}
\le C,
\end{align}
which together with \eqref{3.22} and \eqref{o1} yields that
\begin{align}\label{hh}
\int_0^{T_0}\|\bar x\tilde{\rho}^R_t\|^2_{L^p}dt
&\le C\int_0^{T_0}\big\|\bar x|u^R||\nabla\rho^R|\big\|^2_{L^p(B_R)}dt
\le C\int_0^{T_0}\big\|\bar x^{1-a}u^R\big\|_{L^\infty(B_R)}^2
\big\|\bar x^a\nabla\rho^R\big\|^2_{L^p(B_R)}dt\le C.
\end{align}
With the estimates \eqref{kq1}--\eqref{hh} at hand, we find that the sequence
$(\tilde{\rho}^R,\tilde{u}^R,\tilde{\theta}^R,\tilde{H}^R)$ converges, up to the extraction of subsequences, to some limit $(\rho,u,\theta,H)$ in some weak sense, that is, as $R\rightarrow \infty,$ we have
\begin{align}\label{kq3}
& \tilde{\rho}^R\bar x\rightarrow \rho\bar  x\ \mbox{in} \ C(\overline{B_N}\times [0,T_0]) \mbox{ for any } N>0, \\
& \tilde{\rho}^R\bar x^a\rightharpoonup  \rho \bar x^a\ \mbox{ weakly* in }L^\infty(0,T_0; L^1 \cap H^1 \cap W^{1,q}), \\
& \nabla \tilde{u}^R\rightharpoonup \nabla u,\
\nabla \tilde{\theta}^R\rightharpoonup \nabla \theta,\
\tilde{H}^R \bar x^{\frac{a}{2}}\rightharpoonup H\bar x^{\frac{a}{2}}\
\mbox{ weakly* in }L^\infty(0,T_0; H^1),\\
& \tilde{H}^R\rightharpoonup H
\mbox{ weakly* in }L^\infty(0,T_0; H^2),\\
& \sqrt{\tilde{\rho}^R}\tilde{u}^R\rightharpoonup \sqrt{\rho}u,\
\sqrt{\tilde{\rho}^R}\tilde{\theta}^R\rightharpoonup \sqrt{\rho} \theta,\
\tilde{H}^R_t\rightharpoonup H_t
\ \mbox{ weakly* in }L^\infty(0,T_0; L^2), \\
& \nabla^2\tilde{u}^R\rightharpoonup \nabla^2 u,\
\nabla^2\tilde{\theta}^R\rightharpoonup \nabla^2\theta\
\ \mbox{ weakly in }L^{\frac{q+1}{q}}(0,T_0; L^q)\cap L^2(0,T_0; L^q),\\
& \tilde{H}_t^R\rightharpoonup  H_t,\
\nabla\tilde{H}^R\bar x^{\frac{a}{2}}\rightharpoonup \nabla H\bar x^{\frac{a}{2}},\
\nabla^2 \tilde{H}^R\bar{x}^{\frac{a}{2}}\rightharpoonup \nabla^2 H\bar{x}^{\frac{a}{2}}\ \mbox{ weakly in } L^2(0,T_0; L^2),\\
& \sqrt{\tilde{\rho}^R} \tilde{u}^R_t\rightharpoonup \sqrt{\rho}u_t,\ \sqrt{\tilde{\rho}^R} \tilde{\theta}^R_t\rightharpoonup \sqrt{\rho}\theta_t\ \mbox{ weakly* in } L^\infty(0,T_0; L^2),\\
& \sqrt{\tilde{\rho}^R}\tilde{u}^R_t\rightharpoonup \sqrt{\rho} u_t,\
 \sqrt{\tilde{\rho}^R} \tilde{\theta}^R_t\rightharpoonup \sqrt{\rho} \theta_t\ \mbox{ weakly in } L^2(0,T_0; L^2),
\\
&  \nabla\tilde{u}^R_t\rightharpoonup\nabla u_t,\ \nabla\tilde{\theta}^R_t\rightharpoonup\nabla\theta_t,\
\nabla\tilde{H}^R_t\rightharpoonup\nabla H_t\ \mbox{weakly in}\ L^2(0,T_0; L^2),
\end{align}
with
\begin{align}\label{kq4}
\rho\bar x^a\in L^\infty(0,T_0; L^1), \quad \inf\limits_{0\le t\le T_0}\int_{B_{2N_0}}\rho(x,t)dx\ge \frac14.
\end{align}
Letting $R\rightarrow \infty$, standard arguments together with \eqref{kq3}--\eqref{kq4} show that
$(\rho,u,\theta,H)$ is a strong solution of \eqref{mhd}--\eqref{n4} on $\mathbb{R}^2\times (0,T_0]$ satisfying \eqref{1.10} and \eqref{l1.2}.

\textbf{Step 2. Uniqueness of strong solutions.}
Let $(\rho,u,\theta,H)$ and $(\bar\rho,\bar u,\bar \theta, \bar H)$ be two strong solutions satisfying \eqref{1.10} and \eqref{l1.2} with the same initial data, and denote
\begin{align*}
\Theta\triangleq\rho-\bar\rho,~U\triangleq u-\bar u,~\Psi\triangleq \theta-\bar \theta,~\Phi\triangleq H-\bar H.
\end{align*}

First, subtracting the mass equation \eqref{mhd}$_1$ satisfied by $(\rho,u)$ and $(\bar\rho,\bar u)$ gives that
\begin{align}\label{5.2}
\Theta_t+\bar u\cdot\nabla\Theta+\Theta\divv \bar{u} +\rho \divv U+ U\cdot\nabla \rho= 0.
\end{align}
Multiplying \eqref{5.2} by $2\Theta\bar{x}^{2r}$ for $r\in (1,\tilde a)$ with $\tilde a=\min\{2,a\}$ and integrating by parts over $\mathbb{R}^2$, we deduce from Sobolev's inequality, \eqref{l1.2}, \eqref{3.v2}, and \eqref{3.22} that
\begin{align*}
\frac{d}{dt}\|\Theta\bar{x}^{r}\|_{L^2}^{2}
&\le C \big(\|\bar{u} \bar{x}^{-\frac12}\|_{L^\infty}+\|\nabla\bar{u} \|_{L^\infty}  \big) \|\Theta\bar{x}^{r}\|_{L^2}^{2}+C\|\rho\bar{x}^r\|_{L^\infty}\|\nabla U\|_{L^2}\|\Theta\bar{x}^{r}\|_{L^2}\\
& \quad +C\|\Theta\bar{x}^{r}\|_{L^2}
\|U\bar{x}^{-(\tilde a-r)}\|_{L^{\frac{2q}{(q-2)(\tilde a-r)}}}
\|\bar{x}^{\tilde a} \nabla\rho\|_{L^{\frac{2q}{q-(q-2)(\tilde a-r)}}} \\
&\le C\big(1+\|\nabla\bar u\|_{W^{1,q}}\big)
\|\Theta\bar{x}^{r}\|_{L^2}^{2}+C\|\Theta\bar{x}^r\|_{L^2}
\big(\|\nabla U\|_{L^2}+\|\sqrt{\rho} U\|_{L^2}\big).
\end{align*}
This combined with Gronwall's inequality shows that, for all $0\le t\le T_{0}$,
\begin{align}\label{5.1}
\|\Theta\bar{x}^{r}\|_{L^2}
\le C\int_{0}^{t}\big(\|\nabla U\|_{L^2}+\|\sqrt{\rho} U\|_{L^2}\big)ds.
\end{align}

Next, subtracting \eqref{mhd}$_2$ and \eqref{mhd}$_4$ satisfied by $(\rho,u,\theta,H)$ and $(\bar\rho,\bar u,\bar\theta,\bar H)$ leads to
\begin{align}\label{5.5}
& \rho U_t+\rho u\cdot\nabla U-\mu\Delta U-\nabla\big((\mu+\lambda)\divv U\big) \notag \\
&=
-\rho U\cdot\nabla\bar u-\Theta(\bar u_t+\bar u\cdot\nabla\bar u)-R\nabla(\rho\Psi+\Theta\bar{\theta}) -\frac12\nabla\left(|H|^2-|\bar H|^2\right)+H\cdot\nabla \Phi+\Phi\cdot\nabla \bar H,\\
\label{hlv7.1}
& \Phi_t-\nu\Delta \Phi=H\cdot\nabla U+\Phi\cdot\nabla\bar u-u\cdot\nabla\Phi-U\cdot\nabla\bar H-H\divv U-\Phi\divv \bar{u}.
\end{align}
Multiplying \eqref{5.5} by $U$ and \eqref{hlv7.1} by $\Phi$, respectively, and adding the resulting equations together, we obtain after integration by parts that
\begin{align}\label{5.6}
& \frac{d}{dt}\int \big(\rho |U|^2+|\Phi|^2\big)dx
+2\int\big(\mu|\nabla U|^2+\nu|\nabla\Phi|^2\big) dx \notag \\
& \le C\big(\|\nabla \bar{u}\|_{L^\infty}+\|\nabla u\|_{L^\infty}\big)
\int\big(\rho |U|^2+|\Phi|^2\big)dx \notag \\
& \quad +C\int|\Theta||U|\big(|\bar u_{t}|+|\bar u ||\nabla\bar{u}|\big)dx +C\|\rho\Psi+\Theta\bar{\theta}\|_{L^2}\|\divv U\|_{L^2}\notag\\
&\quad+\frac{1}{2}\int \big(|H|^2-|\bar{H}|^2\big)\divv Udx-\int\Phi\cdot\nabla U\cdot\bar Hdx-\int H\cdot\Phi\divv Udx-\int U\cdot\nabla \bar H\cdot\Phi dx
\notag \\
& \triangleq C\big(\|\nabla \bar{u}\|_{L^\infty}+\|\nabla u\|_{L^\infty}\big)
\int\big(\rho |U|^2+|\Phi|^2\big)dx+\sum_{i=1}^6 K_i.
\end{align}
By H\"older's inequality,\eqref{1.10}, \eqref{l1.2}, \eqref{3.i2}, \eqref{o1}, and \eqref{5.1}, we get that, for $r\in (1,\tilde a)$,
\begin{align*}
K_1 & \le C\|\Theta \bar{x}^{r}\|_{L^2}\|U\bar{x}^{-\frac{r}{2}}\|_{L^4}
\big(\|\bar u_{t}\bar{x}^{-\frac{r}{2}}\|_{L^4}+\|\nabla\bar{u}\|_{L^\infty}
\|\bar u\bar{x}^{-\frac{r}{2}}\|_{L^4}\big) \notag \\
& \le C(\varepsilon')\big(\|\sqrt{\bar{\rho}}\bar{u}_{t}\|_{L^2}^{2}+(1+\|\bar{\rho}\|_{L^\infty}^2) \|\nabla\bar{u}_{t}\|_{L^2}^{2}+\|\nabla\bar{u}\|_{L^\infty}^{2}\big)
\|\Theta \bar{x}^{r}\|_{L^2}^2 \notag\\
& \quad +\varepsilon'\big(\|\sqrt{\rho} U\|_{L^2}^{2}+\big(1+\|\rho\|_{L^\infty}^2) \|\nabla U\|_{L^2}^2\big)\notag \\ & \le C(\varepsilon)\big(1+\|\nabla  {\bar{u}_t}\|_{L^2}^{2}
+\|\nabla^{2}\bar{u}\|_{L^q}^{2}\big)
\int_{0}^{t}\big(\|\nabla U\|_{L^2}^2+\|\sqrt{\rho}U\|_{L^2}^{2}\big)ds
+\varepsilon\big(\|\sqrt{\rho}U\|_{L^2}^{2}+\|\nabla U\|_{L^2}^2\big).
\end{align*}
It deduces from \eqref{1.10} and \eqref{3.3.1} that
\begin{align*}
K_2 &\le C \|\nabla U\|_{L^2}(\|\rho \Psi\|_{L^2}+\|\theta\bar{x}^{-r}\|_{L^\infty}\|\Theta\bar{x}^{r}\|_{L^2})\notag\\
&\le \varepsilon\|\nabla U\|_{L^2}^2+C(\varepsilon)\|\sqrt{\rho} \Psi\|_{L^2}^2+C(\varepsilon)\int_0^t(\|\nabla U\|_{L^2}^2+\|\sqrt{\rho} U\|_{L^2}^2)ds.
\end{align*}
For the term $K_3$, we have
\begin{align*}
K_3&=\frac{1}{2}\int(H\cdot\Phi+\bar{H}\cdot\Phi)\divv Udx\notag\\
&\le C\big(\|H\|_{L^4}+\|\bar{H}\|_{L^4}\big)\|\Phi\|_{L^4}\|\nabla U\|_{L^2}\notag\\
&\le \varepsilon\|\nabla U\|_{L^2}^2+ \varepsilon\|\nabla \Phi\|_{L^2}^2+C(\varepsilon)\|\Phi\|_{L^2}^2,
\end{align*}
while, for the term $K_4+K_5$, we derive from Gagliardo-Nirenberg inequality that
\begin{align*}
K_4+K_5 \le C \big(\| H\|_{L^4}+\|\bar H\|_{L^4}\big) \|\Phi\|_{L^4}\|\nabla U\|_{L^2}
\le \varepsilon\|\nabla U\|_{L^2}^2+\varepsilon\|\nabla \Phi\|_{L^2}^2+C(\varepsilon)\|\Phi\|_{L^2}^2.
\end{align*}
The last term $K_6$ can be estimated as follows
\begin{align*}
K_6\le & C\|U\bar x^{-a}\|_{L^4}\||\nabla \bar H|^{\frac12}\bar x^{a}\|_{L^4}\||\nabla \bar H|^{\frac12}\|_{L^4}\|\Phi\|_{L^4} \notag \\
\le & C\left(\|\sqrt{\rho}U\|_{L^2}+\big(1+\|\rho\|_{L^\infty}\big)\|\nabla U\|_{L^2}\right)\|\nabla \bar H\bar x^{\frac{a}{2}}\|_{L^2}^{\frac12}\|\Phi\|_{L^4} \notag \\
\le & \varepsilon\left(\|\sqrt{\rho}U\|_{L^2}^2+\|\nabla U\|_{L^2}^2\right)+C(\varepsilon)\|\Phi\|_{L^4}^2 \notag \\
\le & \varepsilon\left(\|\sqrt{\rho}U\|_{L^2}^2+\|\nabla U\|_{L^2}^2\right)+\varepsilon\|\nabla \Phi\|_{L^2}^2+C(\varepsilon)\|\Phi\|_{L^2}^2.
\end{align*}
Inserting the above estimates on $K_i$ into \eqref{5.6} and choosing $\varepsilon$ suitably small, we arrive at
\begin{align}\label{zxd}
&\frac{d}{dt}\big(\|\sqrt{\rho}U\|_{L^2}^2+\|\Phi\|_{L^2}^2\big)+\|\nabla U\|_{L^2}^2+\|\nabla\Phi\|_{L^2}^2\notag\\
&\le C\big(1+\|\nabla\bar{u}\|_{L^\infty}+\|\nabla u\|_{L^\infty}+\|\nabla \bar{u}_t\|_{L^2}^2+\|\nabla^2 \bar{u}\|_{L^q}^2\big)\notag\\
&\quad \times\Big(\|\sqrt{\rho}U\|_{L^2}^2+\|\Phi\|_{L^2}^2+\|\sqrt{\rho}\Psi\|_{L^2}^2+\int_0^t\big(\|\nabla U\|_{L^2}^2+\|\sqrt{\rho}U\|_{L^2}^2\big)ds\Big).
\end{align}

Finally, owing to \eqref{l1.2}, \eqref{3.i2}, and \eqref{o1}, we infer from \eqref{mhd}$_3$ that
\begin{align}\label{sha2}
c_v\big(\rho\Psi_t+\rho u\cdot \nabla\Psi\big)-\kappa\Delta\Psi
&=-c_v(\rho U\cdot\nabla \bar{\theta})-\Theta\big(\bar{\theta}_t+\bar{u}\cdot \nabla \bar{\theta}\big)-R\rho\theta\divv U-R(\rho\Psi+\Theta\bar{\theta})\divv\bar{u}\notag\\
&\quad+\frac{\mu}{2}\big(\nabla U+(\nabla U)^{tr}\big):\big(\nabla u+(\nabla u)^{tr}+\nabla \bar{u}+(\nabla \bar{u})^{tr}\big)\notag\\
&\quad+\lambda\divv U\divv(u+\bar{u})+\nu\curl\Phi\cdot\curl(H+\bar{H}),
\end{align}
which multiplied by $\Psi$ and integration by parts leads to
\begin{align}\label{sha3}
&\frac{d}{dt}\|\sqrt{\rho}\Psi\|_{L^2}^2+\|\nabla\Psi\|_{L^2}^2\notag\\
&\le C\|\nabla \bar{\theta}\|_{L^\infty}\|\sqrt{\rho}U\|_{L^2}\|\sqrt{\rho}\Psi\|_{L^2} + C\|\Theta(|\bar{\theta}_t|+|\bar{u}\cdot\nabla\bar{\theta}|)\Psi\|_{L^1}+C\|\rho\theta|\nabla U|\Psi\|_{L^1} \notag\\
&\quad+C\|(\rho \Psi+\Theta\bar{\theta})\divv \bar{u}\Psi\|_{L^1}+C\|\big(\nabla U+(\nabla U)^{tr}\big)\big(\nabla u+(\nabla u)^{tr}+\nabla \bar{u}+(\nabla \bar{u})^{tr}\big)\Psi\|_{L^1}\notag\\
&\quad+C\|\divv U\divv(u+\bar{u})\Psi\|_{L^1}+C\||\curl\Phi||\curl(H+\bar{H})|\Psi\|_{L^1}\notag\\
& \triangleq C\|\nabla \bar{\theta}\|_{L^\infty}\big(\|\sqrt{\rho}U\|_{L^2}^2+\|\sqrt{\rho}\Psi\|_{L^2}^2\big) +\sum_{i=1}^6 T_i.
\end{align}
We need to control the right-hand side terms in \eqref{sha3}. From \eqref{3.i2}, \eqref{3.22}, and \eqref{o1}, it has that
\begin{align*}
T_1&\le C\|\Theta\bar{x}^r\|_{L^2}\|\Psi\bar{x}^{-\frac{r}{2}}\|_{L^4}\big(\|\bar{\theta}_t\bar{x}^{-\frac{r}{2}}\|_{L^4}+\|\nabla \bar{\theta}\|_{L^4}\|\bar{u}\bar{x}^{-\frac{r}{2}}\|_{L^\infty}\big)\notag\\
&\le C\|\Theta\bar{x}^r\|_{L^2}\big(\|\sqrt{\rho}\Psi\|_{L^2}+(1+\|\rho\|_{L^\infty})\|\nabla \Psi\|_{L^2}\big)\big(1+\|\nabla \bar{\theta}_t\|_{L^2}\big)\notag\\
&\le\tilde{\varepsilon}\|\sqrt{\rho}\Psi\|_{L^2}^2+\tilde{\varepsilon}\|\nabla\Psi\|_{L^2}^2+C\big(1+\|\nabla\bar{\theta}_t\|_{L^2}^2\big)\int_0^t\big(
\|\nabla U\|_{L^2}^2+\|\sqrt{\rho}U\|_{L^2}^2\big)ds.
\end{align*}
 By \eqref{1.10} and \eqref{3.3.2}, it yields that
\begin{align*}
T_2\le C\|\nabla U\|_{L^2}\|\sqrt{\rho}\theta\|_{L^\infty}\|\sqrt{\rho}\Psi\|_{L^2}
\le \tilde{\varepsilon}\|\nabla U\|_{L^2}^2+C\|\sqrt{\rho}\Psi\|_{L^2}^2,
\end{align*}
and
\begin{align*}
T_3&\le C\|\nabla\bar{u}\|_{L^\infty}\|\sqrt{\rho}\Psi\|_{L^2}^2+C\|\Theta \bar{x}^r\|_{L^2}\|\Psi\bar{x}^{-\frac{r}{2}}\|_{L^4}\|\nabla \bar{u}\|_{L^4}\|\bar{\theta}\bar{x}^{-\frac{r}{2}}\|_{L^\infty}\notag\\
&\le C\|\nabla\bar{u}\|_{L^\infty}\|\sqrt{\rho}\Psi\|_{L^2}^2+C\|\Theta \bar{x}^r\|_{L^2}\big(\|\sqrt{\rho}\Psi\|_{L^2}+(1+\|\rho\|_{L^\infty})\|\nabla\Psi\|_{L^2}\big)\notag\\
&\le C\|\nabla\bar{u}\|_{L^\infty}\|\sqrt{\rho}\Psi\|_{L^2}^2+\tilde{\varepsilon}\big(\|\sqrt{\rho}\Psi\|_{L^2}^2+\|\nabla\Psi\|_{L^2}^2\big)
+C(\tilde{\varepsilon})
\int_0^t\big(\|\nabla U\|_{L^2}^2+\|\sqrt{\rho}U\|_{L^2}^2\big)ds.
\end{align*}
Remember Lemma \ref{lemma2.5} and \eqref{1.10}, then the same method as \eqref{sha7} runs that
\begin{align*}
\sum_{i=4}^6 T_i&\le C\big(\|\nabla U\|_{L^2}+\|\nabla \Phi\|_{L^2}\big)\|\Psi \bar{x}^{-\frac b4}\|_{L^8}\notag\\
&\quad\times\Big(\big(\|\nabla u\|_{L^4}+\|\nabla H\|_{L^4}\big)^{\frac{1}{2}}\|\sqrt{Z}\bar{x}^{\frac{b}{2}}\|_{L^2}^{\frac{1}{2}}
+\big(\| \nabla\bar{u}\|_{L^4}+\|\nabla \bar{H}\|_{L^4}\big)^{\frac{1}{2}}
\|\sqrt{\bar{Z}}\bar{x}^{\frac{b}{2}}\|_{L^2}^{\frac{1}{2}}\Big)\notag\\
&\le C\big(\|\nabla U\|_{L^2}+\|\nabla \Phi\|_{L^2}\big)\|\Psi \bar{x}^{-\frac b4}\|_{L^8}\notag\\
&\le C\big(\|\nabla U\|_{L^2}+\|\nabla \Phi\|_{L^2}\big)\big(\|\sqrt{\rho}\Psi\|_{L^2}+(1+\|\rho\|_{L^\infty})\|\nabla\Psi\|_{L^2}\big)\notag\\
&\le\tilde{\varepsilon}\big(\|\sqrt{\rho}\Psi\|_{L^2}^2+\|\nabla\Psi\|_{L^2}^2\big)+C(\tilde{\varepsilon})\big(\|\nabla U\|_{L^2}^2+\|\nabla \Phi\|_{L^2}^2,
\end{align*}
where $\bar{Z}\triangleq\frac{\mu}{2}|\nabla \bar{u}+(\nabla \bar{u})^{tr}|^{2}+\lambda(\divv \bar{u})^2+\nu(\curl \bar{H})^{2}$
and $|\nabla u|+|\divv u|+|\curl H|\le C\sqrt{Z}$.
Inserting the above estimates on $T_i$ into \eqref{sha3} and choosing $\tilde{\varepsilon}$ suitably small, we find that
\begin{align}\label{sha9}
\frac{d}{dt}\|\sqrt{\rho}\Psi\|_{L^2}^2+\|\nabla\Psi\|_{L^2}^2&\le C(\|\nabla U\|_{L^2}^2+\|\nabla \Phi\|_{L^2}^2)+C\big(1+\|\nabla\bar{\theta}_t\|_{L^2}^2+\|\nabla\bar{\theta}\|_{L^\infty}+\|\nabla\bar{u}\|_{L^\infty}\big)\notag\\
&\quad\times\left(\|\sqrt{\rho}U\|_{L^2}^2+\|\sqrt{\rho}\Psi\|_{L^2}^2+\int_{0}^{t}\big(\|\nabla U\|_{L^2}^2+\|\sqrt{\rho}U\|_{L^2}^{2}\big)ds\right).
\end{align}
Denoting
\begin{align*}
G(t)\triangleq \|\sqrt{\rho} U\|_{L^2}^{2}+\|\Phi\|_{L^2}^2+\|\sqrt{\rho}\Psi\|_{L^2}^2
+\int_0^t\big(\|\nabla U\|_{L^2}^2+\|\nabla\Phi\|_{L^2}^2+\|\sqrt{\rho}U\|_{L^2}^2+\|\nabla\Psi\|_{L^2}^2\big)ds.
\end{align*}
Then, multiplying \eqref{zxd} by a large constant and adding it up to \eqref{sha9}, we have
\begin{align*}
G'(t)\le C\big(1+\|\nabla \bar u\|_{L^\infty}+\|\nabla u\|_{L^\infty}+\|\nabla \bar \theta\|_{L^\infty}+\|\nabla \bar \theta_t\|_{L^2}^2+\|\nabla\bar u_t\|_{L^2}^2+\|\nabla^2 \bar{u}\|_{L^q}^2\big)G(t),
\end{align*}
which together with Gronwall's inequality and \eqref{1.10} implies that $G(t)=0$.
This gives that
\begin{align*}
U(x,t)=0,\ \Phi(x,t)=0,\ \Psi(x,t)=0,
\end{align*}
for almost every $(x, t)\in\mathbb{R}^2\times(0, T_0]$. Moreover, one infers from \eqref{5.1} that, for almost every $(x, t)\in\mathbb{R}^2\times(0, T_0]$,
\begin{align*}
\Theta=0.
\end{align*}
Thus we finish the proof of the uniqueness of solutions.
\hfill $\Box$

\end{document}